\documentclass{amsart}
\usepackage{amssymb}
\usepackage{amscd}
\usepackage{verbatim}
\usepackage{epsfig}
\begin{document}
\newcommand{\pa}{\partial}
\newcommand{\CI}{C^\infty}
\newcommand{\dCI}{\dot C^\infty}
\newcommand{\supp}{\operatorname{supp}}
\renewcommand{\Box}{\square}
\newcommand{\ep}{\epsilon}
\newcommand{\Ell}{\operatorname{Ell}}
\newcommand{\WF}{\operatorname{WF}}
\newcommand{\WFb}{\operatorname{WF}_{\bl}}
\newcommand{\diag}{\mathrm{diag}}
\newcommand{\sign}{\operatorname{sign}}
\newcommand{\Ker}{\operatorname{Ker}}
\newcommand{\Ran}{\operatorname{Ran}}
\newcommand{\sH}{\mathsf{H}}
\newcommand{\codim}{\operatorname{codim}}
\newcommand{\Id}{\operatorname{Id}}
\newcommand{\cl}{{\mathrm{cl}}}
\newcommand{\piece}{{\mathrm{piece}}}
\newcommand{\bl}{{\mathrm b}}
\newcommand{\scl}{{\mathrm{sc}}}
\newcommand{\Psib}{\Psi_\bl}
\newcommand{\Psibc}{\Psi_{\mathrm{bc}}}
\newcommand{\Psibcc}{\Psi_{\mathrm{bcc}}}
\newcommand{\Psisc}{\Psi_\scl}
\newcommand{\BB}{\mathbb{B}}
\newcommand{\RR}{\mathbb{R}}
\newcommand{\NN}{\mathbb{N}}
\newcommand{\sphere}{\mathbb{S}}
\newcommand{\codimY}{k}
\newcommand{\dimX}{n}
\newcommand{\cO}{\mathcal O}
\newcommand{\cS}{\mathcal S}
\newcommand{\cF}{\mathcal F}
\newcommand{\cL}{\mathcal L}
\newcommand{\cH}{\mathcal H}
\newcommand{\cG}{\mathcal G}
\newcommand{\cU}{\mathcal U}
\newcommand{\cM}{\mathcal M}
\newcommand{\loc}{{\mathrm{loc}}}
\newcommand{\comp}{{\mathrm{comp}}}
\newcommand{\Tb}{{}^{\bl}T}
\newcommand{\Vf}{\mathcal V}
\newcommand{\Vb}{{\mathcal V}_{\bl}}
\newcommand{\Vsc}{{\mathcal V}_{\scl}}
\newcommand{\etat}{\tilde\eta}
\newcommand{\Hsc}{H_{\scl}}
\newcommand{\ff}{{\mathrm{ff}}}

\setcounter{secnumdepth}{3}
\newtheorem{lemma}{Lemma}[section]
\newtheorem{prop}[lemma]{Proposition}
\newtheorem{thm}[lemma]{Theorem}
\newtheorem{cor}[lemma]{Corollary}
\newtheorem{result}[lemma]{Result}
\newtheorem*{thm*}{Theorem}
\newtheorem*{prop*}{Proposition}
\newtheorem*{cor*}{Corollary}
\newtheorem*{conj*}{Conjecture}
\numberwithin{equation}{section}
\theoremstyle{remark}
\newtheorem{rem}[lemma]{Remark}
\newtheorem*{rem*}{Remark}
\theoremstyle{definition}
\newtheorem{Def}[lemma]{Definition}
\newtheorem*{Def*}{Definition}

\title{The inverse problem for the local geodesic ray transform}
\author[Gunther Uhlmann and Andras Vasy]{Gunther Uhlmann and Andr\'as Vasy}
\date{October 7, 2012}
\address{Department of Mathematics, University of Washington, 
Seattle, WA 98195-4350, and Department of Mathematics,
University of California, Irvine, 340 Rowland
Hall, Irvine CA 92697}
\email{gunther@math.washington.edu}
\address{Department of Mathematics, Stanford University, Stanford, CA
94305-2125, U.S.A.}
\email{andras@math.stanford.edu}
\thanks{The authors were partially supported by the National Science Foundation under
grant CMG-1025259 (G.U.\ and A.V.) and DMS-0758357 (G.U.) and
DMS-1068742 (A.V.).}

\begin{abstract}
Under a convexity assumption on the boundary we solve a local
inverse problem, namely we show that the geodesic X-ray transform can
be inverted locally in a stable manner; one even has a reconstruction
formula.
We also show that under an assumption on the
existence of a global foliation by strictly convex hypersurfaces the
geodesic X-ray transform is globally injective. In addition we prove stability estimates and propose a layer stripping type algorithm for reconstruction.
\end{abstract}

\maketitle

\section{Introduction}
Let $X$ be a strictly convex domain in a Riemannian manifold $(\tilde
X,g)$ of dimension $\geq 3$. In this paper we consider the local inverse problem for the
geodesic X-ray transform. That is, for an open set $O\subset
\overline{X}$, we call geodesic segments $\gamma$ of $g$ which are
contained in $O$ with endpoints at $\pa X$ {\em $O$-local geodesics};
we denote the set of these by $\cM_O$. Thus, $\cM_O$ is an open subset of
the smooth manifold of all geodesics, $\cM$.
We then define
the {\em local geodesic transform} of a function $f$ defined on
$X$ as the collection $(If)(\gamma)$ of integrals of $f$ along
geodesics $\gamma\in \cM_O$, i.e.\ as the restriction of the X-ray
transform to $\cM_O$.

In order to state our main theorem in concrete terms, it is useful to
introduce some notation. Let $\rho\in\CI(\tilde X)$ be a defining function of $\pa X$,
considered a function on $\tilde X$ (so $\rho>0$ in $X$, $<0$ on
$\tilde X\setminus \overline{X}$, vanishes non-degenerately at $\pa
X$). Our main theorem is an invertibility result for the local
geodesic transform on neighborhoods of $p$ in $\overline{X}$ of
the form $\{\tilde x>-c\}$, $c>0$, where $\tilde x$ is a function with
$\tilde x(p)=0$, $d\tilde x(p)=-d\rho(p)$, see Figure~\ref{fig:convex-1} below.

\begin{thm*}
For each $p\in\pa X$, there exists a function $\tilde x\in\CI(\tilde
X)$ vanishing at $p$ and with $d\tilde x(p)=-d\rho(p)$
such that for $c>0$ sufficiently small, and with
$O_p=\{\tilde x>-c\}\cap\overline{X}$,
the local geodesic transform is injective on
$H^s(O_p)$, $s\geq 0$.

Further, let $H^s(\cM_{O_p})$ denote the restriction of
elements of $H^s(\cM)$ to $\cM_{O_p}$, and
for $\digamma>0$ let
$$
H^s_\digamma(O_p)=e^{\digamma/(\tilde x+c)} H^s=\{f\in
H^s_{\loc}(O_p):\ e^{-\digamma/(\tilde x+c)} f\in
H^s(O_p)\}.
$$
Then for $s\geq 0$
there exists $C>0$ such
that for all $f\in H^s_\digamma(O_p)$,
$$
\|f\|_{H^{s-1}_\digamma(O_p)}\leq C\|If|_{\cM_{O_p}}\|_{H^s(\cM_{O_p})}.
$$
\end{thm*}

\begin{rem*}
Here the constant $C$ is uniform in $c$ for small $c$, and indeed if
we consider the regions $\{\rho\geq\rho_0\}\cap\{\tilde x>-c\}$ with
$|\rho_0|$ and $|c|$ sufficiently small and such that this
intersection is non-empty, the estimate is uniform in both $c$ and
$\rho_0$.

Further, the estimate is also stable under sufficiently small
perturbations of the metric $g$, i.e.\ the constant is
uniform. (Notice that the hypotheses of the theorem are satisfied for
small perturbations of $g$!)
\end{rem*}

We remark that for this result we only need to assume convexity near the point $p.$
This local result is new even in the case that the metric is conformal
to the Euclidean metric. We also point out that we also get a reconstruction method in the form of a Neumann series. See Section~\ref{sec:sc-calc} for more details.

While this large weight $e^{\digamma/(\tilde x+c)}$ means that the
control over $f$ in terms of $If$ is weak at $\tilde x=-c$, the control is uniform in compact
subsets of $O_p$: these weights are
bounded below on $O_p$ by a positive constant, and bounded above on
compact subsets of $O_p$
(in particular at parts of $\pa X$). Here $\digamma>0$ can be taken
small, but not vanishing. Further, $\tilde x$, whose existence is
guaranteed by the theorem,
is such that $\tilde x=-c$ is concave
from the side of $O_p$.

As an application, we consider domains with compact closure $\overline{X}$ equipped with a function
$\rho:\overline{X}\to[0,\infty)$ whose level sets
$\Sigma_t=\rho^{-1}(t)$, $t<T$, are strictly
convex (viewed from $\rho^{-1}((t,\infty))$ (and $d\rho$ is
non-zero on these level sets), with $\Sigma_0=\pa X$
and $X\setminus\cup_{t\in[0,T)}\Sigma_t=\rho^{-1}([T,\infty))$ either having $0$ measure or having empty interior. (Note in
particular that $\rho$ is a boundary defining function.)

\begin{cor*}
For $X$ and $\rho$ as above, if $X\setminus\cup_{t\in[0,T)}\Sigma_t$
has $0$ measure, the global geodesic transform is
injective on $L^2(X)$, while if it has empty interior, the global geodesic transform is
injective on $H^s(X)$, $s>n/2$.
\end{cor*}

This corollary is an immediate consequence of our main theorem. Indeed,
if $If=0$ and $f\in H^s$, $s>n/2$, $f\neq 0$, then $\supp f$ has non-empty
interior since $f$ is continuous by the Sobolev embedding, while if $f\in L^2$, $f\neq 0$, then $\supp f$ has non-zero
measure.
On the other hand, let $\tau=\inf_{\supp f}\rho$; if $\tau\leq T$ we
are done, for then $\supp f\subset
X\setminus\cup_{t\in[0,T)}\Sigma_t$. Thus, suppose $\tau>T$, so
$f\equiv 0$ on $\Sigma_t$ for $t<\tau$, but there exists
$q\in\Sigma_\tau\cap\supp f$ (since $\supp f$ is closed and
$\overline{X}$ is compact). Now we use
the main theorem on $\rho^{-1}(\tau,\infty)$ to conclude that a
neighborhood of $q$ is disjoint from $\supp f$ to obtain a
contradiction.

In fact, in this global setting we can even take $\tilde x=-\rho$, and
the uniformity of the constants in terms of $c$ and $\rho_0$, as
stated in the remark after the main theorem directly yields that if
$t<T$ then there exists $\delta=\delta_t>0$ such that if $c,\rho_0\in
(t-\delta_t,t+\delta_t)$ then a stability estimate holds (with a
reconstruction method!) for the
region $\rho^{-1}([\rho_0,c))$. Now in general, for $T'<T$, one can
take a finite open cover of $[0,T']$ by such intervals $(t_j',t_j'')$,
$j=1,\ldots,k$ (with, possibly after some reindexing and dropping some
intervals, $t_1'<0$, $t_k''>T'$, $t_j''\in(t_{j+1}',t_{j+1}'')$), and proceed
inductively to recover $f$ on $\cup_{t\in[0,T']}\Sigma_t$ from its X-ray transform, starting with the
outermost region. More precisely, first, using the theorem, one can
recover the restriction of $f$ to
$\rho^{-1}((-\infty,t_1''))$. Then one turns to the next interval,
$(t_2',t_2'')$, 
and notes there is a reconstruction method for the restriction to
$\rho^{-1}((t_2',t_2''))$
of functions $f_2$
supported in $\rho^{-1}((t_2',+\infty))$ (no support condition needed at the other end,
$t_2''$). One applies this to
$f_2=\phi_2 f$, where $\phi_2$ identically $1$ near $\rho^{-1}([t''_1,+\infty))$, supported in $\rho^{-1}((t_2',+\infty))$; since
$f=(1-\phi_2)f+\phi_2 f$, and one has already recovered $(1-\phi_2)f$,
one also knows the X-ray transform of $\phi_2 f$, and thus the theorem
is applicable. One then proceeds inductively, covering
$\rho^{-1}([0,T'])$ in $k$ steps.
This gives a {\em global} stability estimate, and indeed
a reconstruction method doing a reconstruction layer by layer; that
is, we have
(in principle) developed a layer stripping algorithm for this problem.

The geodesic ray transform is closely related to the boundary rigidity
problem of determining a metric on a compact Riemannian manifold from its
boundary distance function. See \cite{SU, I} for recent reviews. The case
considered here is the
linearization of the boundary rigidity problem in a fixed conformal class.
The standard X-ray transform, where one integrates a function along straight
lines, corresponds to the case of the Euclidean metric and is the basis of
medical imaging techniques such as CT and PET. The case of integration along
more general geodesics arises in geophysical imaging in determining the
inner structure of the Earth since the speed of elastic waves generally
increases with depth, thus curving the rays back to the Earth surface. It
also arises in ultrasound imaging, where the Riemannian metric models the
anisotropic index of refraction. Uniqueness and stability
was shown by Mukhometov \cite{Mu} on simple surfaces, and also for more general families of curves
in two dimensions. The case of geodesics was generalized also for simple manifolds to higher dimensions in
\cite{Mu-R}, \cite{Mu}, \cite{BG}. In dimension $n\ge 3$, the paper \cite{FSU} proves injectivity and stability for the X-ray transform integrating over quite a general class of analytic curves with analytic weights, assuming an additional microlocal condition that includes the case of real-analytic metrics for a class of non-simple manifolds. Reconstruction procedures or inversion formulas have not been proven except in a few cases for instance for a class of symmetric spaces,
see \cite{Hel}, and real-analytic curves \cite{FSU}. Our results
generalize support type theorems to the smooth case for the geodesic X-ray transform given in \cite{K} for simple real-analytic metrics.

The global geometric condition that we are imposing is a natural analog of the
condition $\frac{d}{dr} (r/c(r))>0$ proposed by Herglotz \cite{Her} and
Wiechert and Zoeppritz \cite{WZ} for an isotropic radial sound speed $c(r)$.
In this case the geodesic spheres are strictly convex. It is also satisfied
for negatively curved manifolds.
But this condition allows in principle for conjugate points of the metric. In \cite{SU4} one can find a microlocal study of the geodesic X-ray transform with fold caustics.
A similar condition of foliating by convex hypersurfaces was used in \cite{SU5} to 
satisfy the pseudoconvexity condition needed for Carleman estimates.

We also remark that our approach is a completely new one to uniqueness for the global problem for the geodesic ray transform. The only method up to now, except in the real-analytic category \cite{SU}, has been the use of energy type equalities one introduced by Mukhometov \cite{Mu} and developed by several authors which are now called ``Pestov identities". 

The main theorem is proved by considering an operator $A$ which is
essentially a `microlocal normal operator' for the geodesic ray
transform. Let $\rho$ be a boundary defining function of $X$, i.e.\
$\rho>0$ in $X$, $\rho=0$ at $\pa X$, and $d\rho\neq 0$ at $\pa X$; we
assume that in fact $\rho$ is defined on the ambient space $\tilde X$
as above.
First we choose an initial neighborhood
$U$ of $p$ in $\tilde X$ and a function $\tilde x$ defined on it with $\tilde
x(p)=0$, $d\tilde x(p)=-d\rho(p)$, $d\tilde x\neq 0$ on $U$ with convex level sets from the side
of the sublevel sets and
such that $O_c=\{\tilde x>-c\}\cap\{\rho\geq 0\}$ satisfies
$\overline{O_c}\subset U$ is compact. Such a $\tilde x$ exists as can
be seen by slightly modifying $-\rho$, making the level sets slightly
less convex. We define an
operator $L$ which integrates $If$ over a subset of $\cM_{O_c}$ with a
$\CI$ cutoff, and
consider $A=L\circ I$. We
consider this operator as a map between appropriate function spaces on
$O_c$. It turns out that with the subset of geodesics we choose, the exponential conjugate
$A_\digamma$ of $A$ is a pseudodifferential
operator in Melrose's scattering calculus \cite{RBMSpec}. (The
exponential conjugate corresponds to working with exponentially
weighted spaces for $A$.) We show that $A_\digamma$ is a
Fredholm operator, and indeed that it is invertible for $c$ near $0$.

\begin{figure}[ht]\label{fig:convex-1}
\includegraphics[width=80mm]{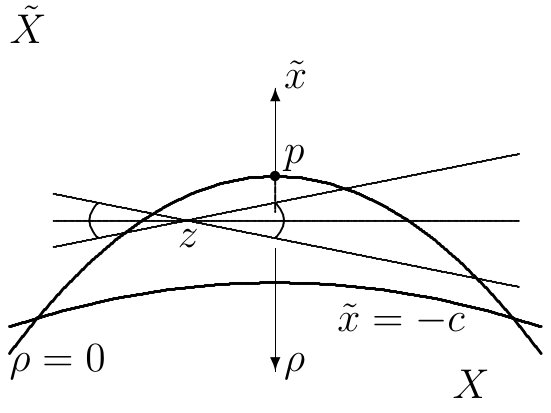}
\caption{The functions $\rho$ and $\tilde x$ when the background is
  flat space $\tilde X$. The intersection of
 $\rho\geq 0$ and $x_c> 0$ (where $x_c=\tilde x+c$, so this is the region $\tilde x>-c$) is the lens
 shaped region $O_p$. Note that, as viewed from the superlevel sets,
 thus from $O_p$, $\tilde
 x$ has concave level sets. At the point $z$, $L$ integrates over
 geodesics in the indicated small angle. As $z$ moves to the
 artificial boundary $x_c=0$, the angle of this cone shrinks like $C
 x_c$ so that in the limit the geodesics taken into account become
 tangent to $x_c=0$.}
\end{figure}

Before giving more details, recall that Stefanov and Uhlmann \cite{SU3} have
shown that under a microlocal condition on the geodesics, one can
recover the singularities of functions from their X-ray transform, and
indeed from a partial X-ray transform (where only some geodesics are
included in the X-ray family $\cM'$). (In
fact, they also showed analogous statements for the transforms on
tensors.) Roughly speaking what one needs is that given a covector
$\nu=(z,\zeta)$, one needs to have a geodesic in $\cM'$ normal to $\zeta$ at $z$
such that in a neighborhood of $\nu$ a simplicity condition is
satisfied. Indeed, under these assumptions, a microlocal version of
the normal operator, $(QI)^*(QI)$, where $Q$ microlocalizes to $\cM'$
roughly speaking, is an elliptic pseudodifferential operator. Now, in
dimension $\geq 3$, if the boundary $\pa X$ is convex,
one can use geodesics which are almost tangent to
$\pa X$ to give a family $\cM'$ which satisfies the above conditions
for $\nu$ with $z$ near $\pa X$. While this gives a recovery of
singularities for the local problem we are considering, it yields no
invertibility or reconstruction. Indeed for the latter we would like
to have $(QI)^*(QI)$ to be an invertible operator on a space of
functions on $O_c$; in particular, as one approaches $\tilde x=-c$ one
would need to only allow integrals over geodesics in a narrow cone,
becoming tangent to $\tilde x=-c$, which takes one outside the
framework of standard pseudodifferential operators.

To remedy this, we introduce the artificial boundary $\tilde x=-c$,
and work with pseudodifferential operators in $x_c=\tilde x+c>0$ which
degenerate at $x_c=0$. Suppressing the $c$ dependence of $x$, the
particular degeneration we end up with is Melrose's scattering
calculus as already mentioned. This is defined on manifolds with
boundary, with boundary defining function $x$, and is based on
degenerate vector fields $x^2\pa_x$ and $x\pa_{y_j}$, where the $(x,y_1,\ldots,y_{n-1})$
are local coordinates. This has the effect of pushing $x=0$ `to
infinity' (these vector fields are complete under the exponential
map). Thus, ultimately, our approach is based on working in a
framework with an artificial boundary which is effectively `at
infinity', and we work with function spaces allowing exponential
growth at this boundary. Thus the control at $x=0$ will be quite weak
in a sense, though one has the standard control when $x$ is bounded
away from $0$. Since $x=0$ is just an artificial boundary, this is
a satisfactory situation.

In fact, for most of the paper we work in a much more general
setting. We consider a family of curves $\gamma_\nu:I\to \tilde X$ parameterized by
$\nu=(z,\zeta)\in S \tilde X$ (the sphere bundle of $\tilde X$ realized as a subbundle
of $T\tilde X$, e.g.\ via a Riemannian metric) with $\gamma_\nu'(0)=\nu$ and
we assume that if $\nu$ is tangent to a level set of $\tilde x$ in
$O_c$, i.e.\ if $\frac{d}{dt}(\tilde x\circ\gamma_\nu)|_{t=0}=0$,
then $\frac{d^2}{dt^2}(\tilde x\circ\gamma_\nu)|_{t=0}\geq C>0$. By
possibly shrinking $U$, we may always assume this in our setting; the
lower bound on the second derivative is a concavity statement for the level sets of
$\tilde x$
from the side of the superlevel sets. Let $x=x_c=\tilde x+c$ as above. Thus,
$x$ is a boundary defining function for $\{\tilde x>-c\}$; for the
time being we regard $c$ as fixed. A
consequence of our uniform concavity statement is that, with
$\lambda=\frac{d}{dt}(x\circ\gamma_\nu)|_{t=0}$, if $C_1>0$ is
sufficiently small and $|\lambda|<C_1
\sqrt{x}$, then $\gamma_\nu$ remains in $x\geq 0$. Rather than using
this range of $\lambda$, we instead use the stronger bound
$|\lambda|< C_2 x$, and define $A$ to be an average:
$$
Af(z)=x^{-1}\int If(\gamma_\nu)\chi(\lambda/x)\,d\mu(\nu),
$$
where $\mu$ is a non-degenerate smooth measure on $S\tilde X$, and $\chi$ has
compact support. We show that for $\digamma\in\RR$,
$$
A_\digamma=x^{-1}e^{-\digamma/x}A e^{\digamma/x}\in\Psisc^{-1,0}(\{x\geq 0\}),
$$
where $\Psisc$ stands for the scattering calculus of Melrose,
and is elliptic in the sense that the standard principal
symbol is such near the boundary (up to the boundary, $x=0$). However, even when this holds
globally on a compact space, this ellipticity is not sufficient
for Fredholm properties (between Sobolev spaces of order shifted by
$1$), or the corresponding estimates, due to
the boundary $x=0$. In general, scattering pseudodifferential
operators also have a principal symbol at the boundary, which is a
(typically non-homogeneous)
function on a cotangent bundle; this needs to be invertible (non-zero) globally
to imply Fredholm properties. Similarly, estimates implying the finite dimensionality
of localized (in $O$) non-trivial nullspace as well as stability estimates,
follow
if this principal symbol is also invertible on $O$. (Note that here
localization {\em does allow} the support in $\{x\geq 0\}$ to include
points at $x=0$!) We thus show that in the case when
$\frac{d^2}{dt^2}(\tilde x\circ\gamma_\nu)|_{t=0}$ is a quadratic form
in $\zeta$ subject to $\lambda=0$, which is the case with geodesics,
for suitable choices of $\chi$, namely essentially cutoff
Gaussians, this principal symbol is invertible when the weight
$\digamma$ satisfies
$\digamma>0$. This implies that
$A_\digamma$ is Fredholm on this space, i.e.\ $A$ itself is Fredholm on
exponentially weighted spaces, where {\em exponential growth} is
allowed at $x=0$. We now recall that $x=x_c$ depends on $c$, with all
estimates uniform for $c$ remaining in a compact set, and the
argument is finished by showing that for $c>0$ sufficiently small one
not only has Fredholm properties but also invertibility, essentially
as the Schwartz kernel has small support.

We note that the geodesic nature of the curves was only used in the
crucial step of showing
that the principal symbol at the boundary is invertible. While our
argument relied on properties of the geodesics to analyze this symbol,
it may well be possible to analyze it in general and prove the result
for more general families of curves. We remark that J. Boman has given in
\cite{Bo}  counterexamples for local uniqueness for the X-ray transform that integrates along lines with a dense family of smooth weights so that we expect some restrictions on the family of curves. 

\section{Scattering calculus}\label{sec:sc-calc}
Melrose's algebra of scattering pseudodifferential operators
$\Psisc^{m,l}(\overline{M})$ on a compact manifold with boundary $\overline{M}$,
see \cite{RBMSpec}, can be thought of either via reducing to a model
on $\RR^n$ (via appropriate charts on $M$, the interior of
$\overline{M}$),
or via a geometric definition. Both are of use in the
current paper; the $\RR^n$ version makes the simplicity of this
algebra transparent, while the geometric definition emphasizes that
infinity in the $\RR^n$-picture is not really `remote', and indeed in
our setting the artificial boundary $\tilde x=-c$ plays $\pa M$, i.e.\
infinity is at a decidedly finite place (moving it to infinity is what
is artificial).

First we start with the $\RR^n$ picture, which is
straightforward. Indeed, the scattering algebra in this setting is a
special case of H\"ormander's Weyl calculus \cite[Section~18.5]{Hor},
which in this particular case has also been studied by Parenti
\cite{Pa} and
Shubin \cite{Sh}. That is, scattering symbols of order $(m,l)$ are defined to be
functions on $\RR^n_z\times\RR^n_\zeta$ satisfying
$$
|D^\alpha_z D^\beta_\zeta a(z,\zeta)|\leq C_{\alpha\beta}\langle z\rangle^{l-|\alpha|}\langle\zeta\rangle^{m-|\beta|},
$$
i.e.\ they are `product type' symbols in $z$ and $\zeta$. Note that
our order convention for the second order $l$, indicating growth/decay in
$z$, is the opposite of that of Melrose \cite{RBMSpec} (i.e.\ our $l$
is $-l$ in \cite{RBMSpec}); we make this deviation so that the
symbol class increases both with $m$ and $l$, i.e.\ so that the two
indices play a parallel role. Their set is
denoted by $S^{m,l}(\RR^n,\RR^n)$ or simply $S^{m,l}$. One then
defines $\Psisc^{m,l}(\RR^n)$ to consist of, say, left quantizations of such
symbols, i.e.\ of operators of the form
\begin{equation}\label{eq:left-quantization}
Au(z)=(2\pi)^{-n}\int e^{i(z-z')\cdot\zeta} a(z,\zeta) u(z')\,dz'\,d\zeta,
\end{equation}
understood as an oscillatory integral.
Right quantizations could be used equally well, i.e.\ one gets the
same class of operators if $a\in S^{m,l}$ but one substitutes
$a(z',\zeta)$ into the oscillatory integral in place of $a(z,\zeta)$. Note that for
$l\leq 0$, $\Psisc^{m,l}(\RR^n)$
is a subspace of H\"ormander's uniform algebra $\Psi_\infty^m(\RR^n)$,
i.e.\ where the above estimates hold without the factor $\langle
z\rangle^{l-|\alpha|}$, and the general weight barely affects the
standard arguments with pseudodifferential operators. The space
$\Psisc^{*,*}(\RR^n)$ is a filtered *-algebra under composition of
operators and taking adjoints (relative to the Euclidean metric), i.e.
$$
A\in\Psisc^{m,l}(\RR^n),\ B\in\Psisc^{m',l'}(\RR^n)\Rightarrow AB\in\Psisc^{m+m',l+l'}(\RR^n)
$$
and
$$
A\in\Psisc^{m,l}(\RR^n)\Rightarrow A^*\in\Psisc^{m,l}(\RR^n).
$$
Further, we define the principal symbol of $A$ to be the equivalence
class of the amplitude
$a$ in \eqref{eq:left-quantization} in $S^{m,l}/S^{m-1,l-1}$, which thus
captures $A$ modulo $\Psisc^{m-1,l-1}(\RR^n)$, i.e.\ one order lower
operators both in terms of the differential order and growth at infinity. With this
definition, the principal symbol of $AB$ is the product of the
principal symbols of $A$ and $B$, while that of $A^*$ is the complex
conjugate of the principal symbol of $A$. In particular, if $A$ is elliptic,
i.e.\ its principal symbol is invertible in the sense that there is
$b\in S^{-m,-l}$ such that $ab-1\in S^{-1,-1}$ (which is independent
of the choice of representative for the principal symbol), then the
standard parametrix construction produces $B\in\Psisc^{-m,-l}(\RR^n)$ such
that $AB-\Id\in\Psisc^{-\infty,-\infty}(\RR^n)$. Operators $R$ in
$\Psisc^{-\infty,-\infty}(\RR^n)$ have a Schwartz function on
$\RR^{2n}$ for their Schwartz kernel; this is just the inverse Fourier
transform of their amplitude $r$ in the $\zeta$ variable evaluated at
$z-z'$ (where $\RR^{2n}=\RR^n_z\times\RR^n_{z'}$, with $z$ the left
and $z'$ the right variable). In particular, such operators are
compact between all {\em polynomially weighted Sobolev spaces}
$H^{s,r}=\langle z\rangle^{-r}H^s(\RR^n)$. Further,
$A\in\Psisc^{m,l}(\RR^n)$ is bounded $H^{s,r}\to H^{s-m,r-l}$, and
if $A$ is elliptic then the parametrix construction and the
compactness we observed shows that $A$ is Fredholm --  it has closed
range, finite dimensional kernel and cokernel, and corresponding
estimates,
$$
\|u\|_{H^{s,r}}\leq
C(\|Au\|_{H^{s-m,r-l}}+\|Fu\|_{H^{-N,-N}}),
$$
where $F$ can be taken a
finite rank
element of $\Psisc^{-\infty,-\infty}(\RR^n)$, and $N$ can be taken
arbitrary.

In order to relate $\Psisc(\RR^n)$ to the geometric setting, and also
in order to explain its classical subalgebra, it is useful to {\em
  compactify} $\RR^n$. Concretely, we compactify $\RR^n$ to a closed
ball $\overline{\RR^n}$ by adding the sphere at infinity $\sphere^{n-1}$. Thus,
$\RR^n\setminus\{0\}$ can be identified with
$(0,\infty)_r\times\sphere^{n-1}_\theta$ via `polar coordinates',
$(r,\theta)\mapsto r\theta$; letting $x=r^{-1}$ we have `reciprocal
polar coordinates', $(0,\infty)_x\times\sphere^{n-1}_\theta$ which allow us to glue a sphere to $x=0$
(corresponding to $r=\infty$) by extending the range of $x$ to
$[0,\infty)$. (Thus, formally, $\overline{\RR^n}$ is the disjoint
union of $\RR^n$ with $[0,\infty)\times\sphere^{n-1}$ modulo the
identification of $\RR^n\setminus\{0\}$ with
$(0,\infty)_x\times\sphere^{n-1}_\theta$.)  Notice that $x=r^{-1}$ is
a boundary defining function near $\pa\overline{\RR^n}$; modifying it
near $0$ gives a global boundary defining function $\rho$. It is straightforward to
check that Schwartz functions on $\RR^n$ are exactly the restrictions
to $\RR^n$
of $\CI$ functions on $\overline{\RR^n}$ which vanish with all
derivatives at $\pa \overline{\RR^n}$. Further, writing $z$ as the
variable on $\RR^n$, the linear vector fields $z_j\pa_{z_k}$ on
$\RR^n$ lift (automatically uniquely, as $\RR^n$ is the interior of
$\overline{\RR^n}$)
to smooth vector fields on $\overline{\RR^n}$ which are tangent to the
boundary, and indeed all smooth vector fields tangent to the boundary
are, away from the origin, linear combinations of these lifts with
coefficients that are smooth on $\overline{\RR^n}$. Since being a
symbol on $\RR^n$, i.e.\ satisfying estimates $|D_z^\alpha a(z)|\leq
C_\alpha \langle z\rangle^{l-|\alpha|}$, is equivalent (away from the
origin, near which one has smoothness) to satisfying stable
estimates under linear vector fields, i.e.\ that $|V_1\ldots V_k a
|\leq C\langle z\rangle^{l}$ for all $k$ and linear vector fields
$V_j$ (with $C$ depending on these), it follows that the lift of a
symbol is a conormal function, i.e.\ a function that
satisfies $\rho^{l}V_1\ldots V_k a\in L^\infty$ whenever $V_j$ are
vector fields tangent to $\pa\overline{\RR^n}$, and conversely, every
conormal function is the lift of a symbol. Correspondingly
$\rho^{-l}\CI(\overline{\RR^n})\subset S^l(\RR^n)$; these are the
`classical' or `one-step' symbols; the Taylor series of a $\CI$
function at the boundary gives rise to the expansion (with $x=\rho$
near $x=0$)
$$
\sum_{j\geq 0}x^{-l+j} a_j(\omega)=\sum_{j\geq 0} r^{l-j}a_j(\omega),
$$
understood as an asymptotic sum.

One can now compactify {\em each factor} of $\RR^n_z\times\RR^n_\zeta$ to
define the compactified space of scattering symbols
$\overline{\RR^n}\times\overline{\RR^n}$; we write $\rho_\pa$ for
the boundary defining function in the first factor (`position', $z$) and $\rho_\infty$ for
that in the second factor (`momentum', $\zeta$). The same considerations as
above show that a scattering symbol on $\RR^n\times\RR^n$ of order
$(m,l)$ corresponds to a conormal function on
$\overline{\RR^n}\times\overline{\RR^n}$, i.e.\ one satisfying
$\rho_\infty^m\rho_\pa^{l}V_1\ldots V_k a\in L^\infty$ whenever $V_j$ are
vector fields tangent to both boundary hypersurfaces of
$\overline{\RR^n}\times \overline{\RR^n}$.
Classical symbols, as before, then are elements of
$\rho_\pa^{-l}\rho_\infty^{-m}\CI(\overline{\RR^n}\times\overline{\RR^n})$,
i.e.\ functions of the form $a=\rho_\pa^{-l}\rho_\infty^{-m} \tilde a$,
$\tilde a\in \CI(\overline{\RR^n}\times\overline{\RR^n})$. Note
that for a classical symbol, its equivalence class in
$S^{m,l}/S^{m-1,l-1}$ can be represented by
$\rho_\pa^{-l}\rho_\infty^{-m}$ times the function $a_0=\tilde a|_{\pa(\overline{\RR^n}\times\overline{\RR^n})}$ on
$\pa(\overline{\RR^n}\times\overline{\RR^n})$ in the sense that any smooth
extension $\tilde a'$ of this function to
$\overline{\RR^n}\times\overline{\RR^n}$ produces an element of the
equivalence class of $a$. Ellipticity then simply means the
non-vanishing of this function $a_0$. Note also that this principal
symbol can be thought of as consisting of two parts, namely the
standard principal symbol, at $\rho_\infty=0$, and the
`boundary principal symbol' at $\rho_\pa=0$. We also write
$\Psisc^{m,l}(\overline{\RR^n})=\Psisc(\RR^n)$.

If $\overline{M}$ is a manifold with boundary with interior $M$, we
can now define $\Psisc^{m,l}(\overline{M})$, much as the standard
pseudodifferential algebra is defined on manifolds by locally
identifying the manifold with $\RR^n$ and imposing that on such charts
$U\times U$ the Schwartz kernel of the operator is that of a
pseudodifferential operator on $\RR^n$, and allowing additional globally
smooth terms in the Schwartz kernel. In our case, the analogous
construction is locally identifying
$\overline{M}$ with $\overline{\RR^n}$, and imposing that on such charts
$U\times U$ the Schwartz kernel of the operator is that of an element
of $\Psisc^{m,l}(\overline{\RR^n})$, and allowing additional globally
{\em Schwartz} (i.e.\ rapidly decaying with all derivatives, smooth)
terms in the Schwartz kernel.
As in the standard manifold case, all the basic properties of the
algebra generalize (one needs to impose some proper support conditions
in the absence of compactness). Concretely, the weighted Sobolev
spaces $\Hsc^{s,r}(\overline{M})$ are also defined by local
identification with $\overline{\RR^n}$, and then $A\in
\Psisc^{m,l}(\overline{\RR^n})$ implies that $A$ is bounded from
$\Hsc^{s,r}(\overline{M})$ to $\Hsc^{s-m,r-l}(\overline{M})$.

It is also of some use to work out the behavior of the Schwartz kernel
of elements of $\Psisc^{m,l}(\overline{M})$ on
$\overline{M}\times\overline{M}$. In view of the previous definition,
this reduces to a calculation for $\Psisc^{m,l}(\overline{\RR^n})$
(modulo Schwartz terms which we ignore as they give elements of
$\Psisc^{-\infty,-\infty}(\overline{M})$). Thus, using local
coordinates $y$ on $\sphere^{n-1}$, and corresponding coordinates
$(x,y,x',y')$ on $\overline{M}\times\overline{M}$, one checks that in the coordinates
$$
x,\ y,\ X=\frac{x-x'}{x^2},\ Y=\frac{y-y'}{x},
$$
valid for $x>0$,
so the diagonal is $X=0$, $Y=0$ when $x>0$,
the Schwartz kernel of an element of $\Psisc^{m,l}$ is of the form
$x^{-l}\tilde K$, where $\tilde K$ is smooth in $(x,y)$ down to
$x=0$ with values in conormal distributions on $\RR^n_{X,Y}$, conormal
to $\{X=0,\ Y=0\}$, which are Schwartz at infinity (i.e.\ decay
rapidly at infinity with all derivatives). Further, the boundary
principal symbol is simply $x^{-l}$ times the Fourier transform in $(X,Y)$
of $\tilde K|_{x=0}$ (a restriction which makes sense in view of the
stated smoothness). In particular, when $l=0$, we need to check that
$$
\int e^{-i\xi X-i\eta\cdot Y} \tilde K(0,y,X,Y)\,dX\,dY
$$
is a non-zero function of $(y,\xi,\eta)$, with a lower bound
$C\langle(\xi,\eta)\rangle$, $C>0$, for its absolute value (which
means we also need a uniform bound at infinity in addition to the
invertibility). Checking this will be the main step of the arguments
presented in the next section.

We mention here that vector fields in
$\Vsc(\overline{M})=x\Vb(\overline{M})$, where $\Vb(\overline{M})$ is
the set of all smooth vector fields tangent to $\pa M$, are in
$\Psisc^{1,0}(\overline{M})$, and indeed the Sobolev spaces of
positive integer differential orders $s$ are equivalently defined by
$u\in\Hsc^{s,r}(\overline{M})$ if and only if $x^{-r}V_1\ldots V_k u\in
L^2_{\scl}(\overline{M})$ for all $k\leq s$ (including $k=0$) and
$V_j\in \Vsc(\overline{M})$; here $L^2_{\scl}(\overline{M})$ is the
$L^2$ space given by identification by $\overline{\RR^n}$, i.e.\ the
measure (or density) is, up to a non-degenerate positive multiple,
$r^{n-1}\,dr\,dy=x^{-n-1}\,dx\,dy$. (Densities like this may be called
{\em scattering densities}.)

We now briefly relate the standard Sobolev spaces $H^s(\overline{M})$
to $\Hsc^{s,r}(\overline{M})$ for $s\geq 0$. First, for $s=0$, the
above description gives $\Hsc^{0,-(n+1)/2}(\overline{M})=L^2
(\overline{M})$ (in the sense of equivalent norms). Next, using that
$V'_1\ldots V'_k u\in L^2 (\overline{M})$ for $k\leq s$ and $V'_j$
smooth vector fields on $\overline{M}$ (which is equivalent to $u\in
H^s(\overline{M})$)
implies that $V_1\ldots V_k u\in
\Hsc^{0,-(n+1)/2} (\overline{M})$ for $k\leq s$ and
$V_j\in\Vsc(\overline{M})$ (since all elements of $\Vsc(\overline{M})$
are smooth vector fields), i.e.\ that $u\in
\Hsc^{s,-(n+1)/2}(\overline{M})$, so
\begin{equation}\label{eq:reg-Sob-in-sc}
H^s(\overline{M})\subset \Hsc^{s,r}(\overline{M}),\ r\leq -\frac{n+1}{2},
\end{equation}
with continuous inclusion map. For the converse direction, we note that if
$V'$ is a smooth vector field, then $x^2 V'\in\Vsc(\overline{M})$. Thus,
$V_1\ldots V_k u\in
\Hsc^{0,2s-(n+1)/2} (\overline{M})$ for $k\leq s$ and
$V_j\in\Vsc(\overline{M})$, so $x^{-2s}V_1\ldots V_k u\in
\Hsc^{0,-(n+1)/2} (\overline{M})$, so $x^{-2}V_1\ldots x^{-2}V_k u\in
\Hsc^{0,-(n+1)/2} (\overline{M})$, implies that
$V'_1\ldots V'_k u\in L^2 (\overline{M})$ for $k\leq s$ and $V'_j$
smooth vector fields. Thus,
\begin{equation}\label{eq:sc-Sob-in-reg}
\Hsc^{s,r}(\overline{M})\subset H^s(\overline{M}),\ r\geq -\frac{n+1}{2}+2s,
\end{equation}
with continuous inclusion map.
There are similar inclusions between negative order spaces. For
instance, as $H^{-s}(\overline{M})=(H^{s}_0(\overline{M}))^*$, $s\geq 0$, via
identification by the $L^2$ pairing, and as $H^s_0(\overline{M})$ is a closed
subspace of $H^s(\overline{M})$, the inclusion
\eqref{eq:reg-Sob-in-sc}
gives the continuous inclusion map on the dual spaces
\begin{equation}\label{eq:neg-sc-Son-in-reg}
\Hsc^{-s,-r}(\overline{M})\subset H^{-s}(\overline{M}),\ -r\geq \frac{n+1}{2}.
\end{equation}

Finally we discuss what happens when ellipticity holds only
locally. Thus, suppose $O$ is an open subset of $\overline{M}$ on
which $A\in\Psisc^{m,l}(\overline{M})$ is elliptic, and suppose that
$K\subset O$ is a compact subset. Let $\phi$ be supported in $O$,
identically $1$ on $K$; let $O'$ be a neighborhood of $\supp\phi$ with
closure compactly contained in $O$. By the ellipticity assumption, there is a {\em
  local parametrix} $G\in\Psisc^{-m-,l}(\overline{M})$ for $A$ such
that $GA=\Id+E$, $E\in\Psisc^{0,0}(\overline{M})$, but over $O'$ the
better conclusion that $E$ is, locally, in $\Psisc^{-\infty,-\infty}$,
holds, so $\phi E\phi\in
\Psisc^{-\infty,-\infty}(\overline{M})$. Thus, $\phi E\phi$ is compact
on any polynomially weighted Sobolev space, so in particular there is
a finite rank operator $F\in\Psisc^{-\infty,-\infty}(\overline{M})$
supported in $O\times O$ such that $\Id+\phi E\phi-F$ is invertible.
Now suppose that $v$ is supported in $K$, so $\phi
v=v$. Then $\phi GA\phi=\phi^2+\phi E\phi$ shows that
$$
(\Id+\phi E\phi)v=\phi G Av,
$$
so
$$
v=(\Id+\phi E\phi-F)^{-1}\phi GAv-(\Id+\phi E\phi-F)^{-1}Fv.
$$
In particular, if $Av=0$ then $v$ is in a finite dimensional space,
namely the range of $(\Id+\phi E\phi-F)^{-1}F$, and if one chooses a
complementary subspace $V$ of $\Ker A\cap\{w:\ \supp w\subset K\}$ in a weighted Sobolev space
$\Hsc^{s,r}(\overline{M})$, then there is a constant $C>0$ such that for $v\in V\cap\{w:\ \supp w\subset
K\}$,
$$
\|v\|_{\Hsc^{s,r}(\overline{M})}\leq C\|Av\|_{\Hsc^{s-m,r-l}(\overline{M})},
$$
i.e.\ a stability estimate holds.

Now suppose that one has a family of operators, $A_t$, $t\in[0,T]$,
depending continuously on $t$ in $\Psisc^{m,l}(\overline{M})$, with
each element of the family being elliptic on $O$ (and thus there is a uniform
constant in the estimates over compact subsets of $O$). Suppose also
that we have a continuous function $f$ on $[0,T]$ with $f(0)=0$, a
compact subset $K_0$ of $O$, and
a family of open sets $\cO_t$, $t>0$ in $K_0$, with the
boundary defining function satisfying $x\leq f(t)$ on $\cO_t$, and we are
interested in distributions $v$ supported in $\cO_t$. In view of the
uniform elliptic estimates, choosing $O'$ a neighborhood of $K_0$
with closure compactly contained in $O$, we then have families of operators $G_t$
and $E_t$, depending continuously on $t\in[0,T]$, with values in
$\Psisc^{-m,-l}(\overline{M})$, resp.\ $\Psisc^{0,0}(\overline{M})$,
such that on $O'$, $E_t$ is uniformly in
$\Psisc^{-\infty,-\infty}(\overline{M})$. Thus, the Schwartz kernel
$K_t$ of $E_t$ satisfies that for any $N$, $x^{-N}(x')^{-N}K_t$ is bounded (with
values in scattering densities in the right, i.e.\ primed, factor),
i.e.\ locally is of the form $\kappa_t \frac{dx'\,dy'}{(x')^{n+1}}$
with $|\kappa_t(x,y,x',y')|\leq C_N x^N (x')^N$. (Notice that the fact
that we used `scattering' densities is thus of little relevance; any
polynomial factor such as $(x')^{-n-1}$, can make no difference.) If
$\phi_t\in\CI_c(\overline{M})$ is now supported in $\cO_t$ and takes
values in $[0,1]$, then $\phi_t E_t\phi_t$ has kernel
$\phi_t(x,y)\phi_t(x',y')\kappa_t \frac{dx'\,dy'}{(x')^{n+1}}$, with
$|\phi_t(x,y)\phi_t(x',y')\kappa_t|\leq C'_N f(t)^{2N}x^{n+1}
(x')^{n+1}$ for all $N$, and thus by Schur's lemma is bounded on
$L^2_{\scl}(\overline{M})$ with norm $\leq C''_N f(t)^{2N}$. In
particular, there is $t_0>0$ such that the norm is $<1/2$ for $t\in
(0,t_0]$. Thus, $\Id+\phi_t E_t\phi_t$ is invertible for such $t$, and the
previous arguments give that if $K_t\subset \cO_t$ is compact then for $t\in(0,t_0]$,
$$
\Ker A_t\cap \{w:\ \supp w\subset K_t\}=\{0\}
$$
and for $v$ supported in $K_t$ one has the stability estimate (with
uniform constant $C$)
$$
\|v\|_{\Hsc^{s,r}(\overline{M})}\leq C\|A_t v\|_{\Hsc^{s-m,r-l}(\overline{M})}.
$$
We remark here that $(\Id+\phi_t E_t\phi_t)^{-1}$ can be constructed
by a Neumann series, and thus ultimately our whole argument is
completely constructive.

In our setting we start with an ambient manifold $\tilde X$ with equipped with a function
$\tilde x$ with non-degenerate level sets near the $0$ value, let
$x_c=\tilde x+c$ ($c$ near $0$), let $M_c=\{x_c> 0\}$, identify a
neighborhood of $Y=\{\tilde x=0\}$ with $Y\times
(-\delta,\delta)_{\tilde x}$, and have a
family of operators
$B_c\in\Psisc^{m,l}(\overline{M_c})$ with Schwartz kernel localized in
$\tilde x<c_0$ (in both factors), where $c_0>0$ is small. We further
have a fixed set $O\subset\tilde X$ with compact closure, $K\subset O$
compact, and a
function $f$ continuous on $[0,\delta)$ with $f(0)=0$, such that on
$O\cap M_c$, $x_c\leq f(c)$. In order to analyze the $B_c$ as $c\to 0$,
we regard these instead as
operators on $M_0=\{\tilde x>0\}$ by letting $A_c=(\Phi_c^{-1})^*
B_c\Phi_c^*$, $\Phi_c(\tilde x,y)=(\tilde x+c,y)$ which maps $M_c$ to
$M_0$. The operators $A_c$ obtained by this procedure (with the
parameter being $c$ rather than $t$), together with the corresponding
translates $\cO_c$ and $K_c$ of $O\cap M_c$ and $K\cap M_c$
satisfy all the requirements of
the previous paragraphs, and thus conclusions apply, which, when
translated to $B_c$ give that for sufficiently small $c$
$$
\Ker B_c\cap \{w\in \Hsc^{s,r}(M_c):\ \supp w\subset K\cap M_c\}=\{0\}
$$
and for $v\in\Hsc^{s,r}(M_c)$ supported in $K$ one has the stability estimate (with
uniform constant $C$)
$$
\|v\|_{\Hsc^{s,r}(\overline{M_c})}\leq C\|B_c v\|_{\Hsc^{s-m,r-l}(\overline{M_c})}.
$$

Further, in our setting, the operators $B_c$ are in fact of the form
$$
B_c=x_c^{-1}e^{-\digamma/x_c}L_c
e^{\digamma/x_c}\in\Psisc^{-1,0}(\overline{M_c}),
$$
so we in fact
obtain that for sufficiently small $c$
$$
\Ker L_c\cap \{e^{-\digamma/x_c} w\in \Hsc^{s,r}(M_c):\ \supp w\subset K\cap M_c\}=\{0\}
$$
and for $v\in e^{\digamma/x_c} \Hsc^{s,r}(M_c)$ supported in $K$ one has the stability estimate (with
uniform constant $C$)
$$
\|e^{-\digamma/x_c} v\|_{\Hsc^{s,r}(\overline{M_c})}\leq C\|e^{-\digamma/x_c} L_c v\|_{\Hsc^{s+1,r-1}(\overline{M_c})}.
$$
Notice that this is an exponentially weak estimate at $\pa M_c$, i.e.\
at $x_c=0$, but the exponential factor is immaterial in
$x_c>0$. Notice also that if $\tilde v\in H^s(\tilde X)$, say, then for
$\digamma>0$ its restriction $v$ to $M_c$ is in $e^{\digamma/x_c}
\Hsc^{s,r}(M_c)$ for all $r$, i.e.\ the results are in fact applicable
to $v$.

\section{Proofs}
Suppose first that $X$ is a domain in $(\tilde X,g)$, $p\in\pa X$, and
$\pa X$ is geodesically strictly convex at $p$ (hence near $p$).
That is, with $\rho$ a boundary defining function of $\overline{X}$,
we have (with $G$ the dual metric, and metric function) that for
covectors $\beta\in T^*_p\tilde X\setminus o$,
$$
(H_G\rho)(\beta) =0\Rightarrow (H_G^2\rho)(\beta)<0.
$$
In particular, by compactness of the unit sphere and homogeneity,
there is a neighborhood $U_0$ of $p$ in $\tilde X$ and $C_0>0$, $\delta>0$ such that for
covectors $\beta\in T^*_{U_0}\tilde X\setminus o$,
$$
|(H_G\rho)(\beta)|<\delta G(\beta)^{1/2}\Rightarrow (H_G^2\rho)(\beta)\leq -C_0 G(\beta).
$$
We then want to define a function $\tilde x$ near $p$ such that
$\tilde x(p)=0$, the
region $\tilde x\geq -c$, $\rho\geq 0$, is compact for $c>0$ small,
and the level sets of $\tilde x$ are concave from the side of this
region (i.e.\ the super-level sets of $\tilde x$). By shrinking $U_0$
if needed, we may assume that it is a coordinate neighborhood of $p$. Concretely we let,
for $\ep>0$ to be decided, an with $|.|$ the Euclidean norm,
$$
\tilde x(z)=-\rho(z)-\ep|z-p|^2;
$$
then $\tilde x\geq -c$ gives $\rho+\ep|z-p|^2\leq c$ and thus
$\rho\leq c$; further, with $\rho\geq 0$ this gives $|z-p|^2\leq
c/\ep$. Thus, for $c/\ep$ sufficiently small, the region $\tilde x\geq
-c$, $\rho\geq 0$, is compactly contained in $U_0$. Further, for
$\beta\in T^*_{U_0} \tilde X$,
$H_G\tilde x(\beta)=-H_G\rho(\beta)-\ep H_G|.-p|^2$, so $H_G\tilde
x=0$ implies $|H_G\rho|<C'\ep G^{1/2}$, so with $\delta>0$ as above there is $\ep'>0$ such that
for $\ep\in(0,\ep')$, $H_G\tilde x=0$ in $U_0$ implies $|H_G
\rho|<\delta G^{1/2}$, and then, for $\ep<\ep'$,
$$
H_G^2\tilde
x=-H_G^2\rho-\ep H_G^2|.-p|^2\geq (C_0-C''\ep) G.
$$
Thus, there is $\ep_0>0$ such that for $\ep\in(0,\ep_0)$,
$H_G^2\tilde x\geq (C_0/2) G$ at $T^*_p
\tilde X$ when $H_G\tilde x$ vanishes. Thus taking $c_0>0$ sufficiently
small (corresponding to $\ep_0$), we have constructed a
function $\tilde x$ defined on a neighborhood $U_0$ of $p$ with
concave level sets (from the side of the super-level sets) and such
that for $0\leq c\leq c_0$,
$$
O_c=\{\tilde x>-c\}\cap\{\rho\geq 0\}
$$
has compact closure in $U_0\cap\overline{X}$.

From now on we work with $x_c=\tilde x+c$, which is the boundary
defining function of the region $x_c\geq 0$; we suppress the $c$
dependence and simply write $x$ in place of $x_c$. For most of the
following discussion we completely ignore the actual boundary,
$\rho=0$; this will only play a role at the end since ellipticity
properties only hold in $U_0$ and we need $f$ to be supported in
$\rho\geq 0$, ensuring localization, in order to obtain injectivity
and stability estimates. Thus, completing $\tilde x$ to a coordinate
system $(\tilde x,y)$ on a neighborhood $U_1\subset U_0$ of $p$, for
each point $(\tilde x,y)$ we can parameterize geodesics through this
point by the unit sphere; the relevant ones for us are `almost
tangent' to level sets of $\tilde x$, i.e.\ we are interested in ones
with tangent vector $c(\lambda\pa_x+\omega\pa_y)$, $c>0$ (to say have
unit length), $\omega\in\sphere^{n-2}$, and $\lambda$ relatively
small.

Now, the geodesic corresponding to $(z_0,\nu_0)$, $\gamma=\gamma_{z_0,\nu_0}$, is the
projection of the bicharacteristic $\tilde\gamma$ emanating from
$(z_0,g_{z_0}(\nu_0))=(z_0,\zeta_0)$ (i.e.\ the integral curve of
$H_G$ through this point; here we are using the metric $g_{z_0}$ to turn the
vector $\nu_0$ into a covector)
which thus satisfies
$(\frac{d}{dt}\tilde\gamma)(t)=H_G(\tilde\gamma(t))$, so
$\frac{d^2}{dt^2}(f\circ\tilde\gamma)(t)=H_G^2f(\tilde\gamma(t))$. Thus,
if $f$ is a function on the base space $\tilde X$ then
$(\frac{d^2}{dt^2}\gamma)(0)=(H_G^2
f)(\gamma(0),g_{\gamma(0)}(\gamma'(0)))$. But $H_G^2$ is homogeneous
degree two in the fiber (second) variable of the cotangent bundle, and
it is a polynomial, which shows that
$\frac{d^2}{dt^2}(f\circ\tilde\gamma)(0)$ is a quadratic polynomial in
$\nu$.

We now make this more concrete. For this, we use a fibration by level sets of a function
$x$ with non-vanishing differential. Letting $V$ be a vector field orthogonal with respect to $g$ to
these level sets with $Vx=1$, and using $\{x=0\}$ as the initial
hypersurface, the flow of $V$ (locally) identifies a neighborhood of
$\{x=0\}$ with $(-\ep,\ep)_x\times\{x=0\}$, with the first coordinate
being exactly the function $x$ (since time $t$ flow by $V$ changes the
value of $x$ by $t$). In particular, choosing coordinates $y_j$ on
$\{x=0\}$, we obtain coordinates on this neighborhood such that
$\pa_{y_j}$ and $\pa_x$ are orthogonal, i.e.\ the metric is of the
form $f(x,y)\,dx^2+h(x,y,dy)$, and the dual metric is of the form
$$
F(x,y)\xi^2+\sum H_{ij}(x,y)\eta_i\eta_j,
$$
with $f,F>0$,
so (with $h_{ij}$ denoting the metric components, so $H_{ij}$ is the
inverse matrix of $h_{ij}$),
$$
\frac{dx}{dt}=2F(x,y)\xi,\ \frac{dy_i}{dt}=2\sum H_{ij}(x,y)\eta_j,\
-\frac{d\xi}{dt}=\frac{\pa F}{\pa x}\xi^2+\sum\frac{\pa H_{ij}}{\pa x}(x,y)\eta_i\eta_j,
$$
and thus
\begin{equation*}\begin{aligned}
\frac{1}{2}\frac{d^2 x}{dt^2}
&=2\frac{\pa F}{\pa x}(x,y)F(x,y)\xi^2 +2\sum\frac{\pa F}{\pa
  y_i}H_{ij}(x,y)\eta_j\xi\\
&\qquad-F(x,y)\frac{\pa F}{\pa x}\xi^2
-F(x,y)\sum\frac{\pa H_{ij}}{\pa
  x}(x,y)\eta_i\eta_j
\end{aligned}\end{equation*}
which at $\frac{dx}{dt}=0$, thus $\xi=0$, simplifies to
\begin{equation*}\begin{aligned}
&-\sum\frac{\pa H_{ij}}{\pa
  x}(x,y)h_{ik}(x,y)h_{jl}(x,y)\frac{dy_k}{dt}\frac{dy_l}{dt}\\
&=-\sum\frac{\pa H_{ij}}{\pa x}(x,y)h_{ik}(x,y)h_{jl}(x,y)\omega_k\omega_l.
\end{aligned}\end{equation*}
Here we used the unit sphere for the $\omega$-parameterization. Note
that
$$
-(\sum\frac{\pa H_{ij}}{\pa x}(x,y)h_{ik}(x,y)h_{jl}(x,y))_{kl}
$$ 
is
positive definite by our assumptions. Thus, for geodesics we have a
positive definite quadratic form
$$
\alpha(x,y,\omega,0,0)=
-\sum(\sum\frac{\pa H_{ij}}{\pa x}(x,y)h_{ik}(x,y)h_{jl}(x,y))_{kl}\omega_k\omega_l.
$$

In fact, as explained in the introduction, we mostly work in the
following more general setting.
We consider integrals along a family of $\CI$ curves
$\gamma_{x,y,\lambda,\omega}$ in $\RR^n$,
$(x,y,\lambda,\omega)\in\RR\times\RR^{n-1}\times\RR\times\sphere^{n-2}$,
depending smoothly ($\CI$) on the parameters,
typically (but not necessarily) geodesics.
Here $\RR^{n-1}_y$ could be replaced by
an arbitrary manifold and below we make $x$ small, so effectively we
are working in a tubular neighborhood of a codimension one
submanifold of an arbitrary manifold, such as $\tilde X$. However,
since the changes in the manifold setting are essentially just
notational, for the sake of clarity we work with $\RR^n$. Further,
below we work with neighborhoods of a compact subset $\{0\}\times
K\subset \RR_x\times\RR^{n-1}_y$; $\gamma_{x,y,\lambda,\omega}(t)$ would only need to be
defined for $(x,y)$ in a fixed neighborhood $\tilde U$ of $\{0\}\times
K$ and for $|\lambda|<\tilde\delta_0$, and $|t|<\tilde\delta_0$, $\tilde\delta_0>0$ a fixed constant.

The basic feature we need is that
for $x\geq 0$ and for $\lambda$ sufficiently small, depending on $x$,
the curves stay in $[0,\infty)\times\RR^{n-1}$.
Thus, for $x=0$ only the parameter value
$\lambda=0$ is allowed; in our concrete setting $|\lambda|\leq
C_0\sqrt{x}$ works for suitably small $C_0>0$. However, it is convenient to use an even
smaller range of $\lambda$, such as $|\lambda|\leq C_0x$. So concretely
assume that
\begin{equation*}\begin{aligned}
&\gamma_{x,y,\lambda,\omega}(0)=(x,y),\
\gamma'_{x,y,\lambda,\omega}(0)=(\lambda,\omega),\\
&\gamma''_{x,y,\lambda,\omega}(t)=2(\alpha(x,y,\lambda,\omega,t),\beta(x,y,\lambda,\omega,t)),
\end{aligned}\end{equation*}
and
$$
\alpha(0,y,0,\omega,0)\geq 2C>0,
$$
with $\alpha$, $\beta$ smooth. This implies that if
$K\subset\RR^{n-1}$ is compact, then for a sufficiently
small neighborhood $U$ of $\{0\}\times K$ in $\RR^n$ (with compact closure), and for $\lambda$ and $t$ sufficiently
small, say $|\lambda|, |t|<\delta_0$, where $\delta_0>0$, one has
$$
\alpha(x,y,\lambda,\omega,t)\geq C>0.
$$
One may assume that $x<\delta_0$ on $U$.
Thus, writing $\gamma(t)=(x'(t),y'(t))$,
$$
x'=x+\lambda t+t^2\int_0^1(1-s) \alpha(x,y,\lambda,\omega,s)\,ds\geq
x+\lambda t+Ct^2/2,
$$
so if $|t|<\delta_0$, $(x,y)\in U$, $|\lambda|<\delta_0$ then
\begin{equation}\label{eq:x-xp-comp}
x'\geq \frac{C}{2}\Big(t+\frac{\lambda}{C}\Big)^2+\Big(x-\frac{\lambda^2}{2C}\Big).
\end{equation}
Thus, for $|\lambda|\leq \sqrt{2C}\sqrt{x}$ (and $|\lambda|<\delta_0$), $x'\geq 0$, i.e.\ the curves
remain in the half-space $x'\geq 0$ at least for $|t|<\delta_0$. Further, if we fix $x_0>0$, then
$x'\geq x_0$ provided $|t+\frac{\lambda}{C}|>\sqrt{2x_0/C}$ and
$|t|<\delta_0$, thus when $|\lambda|\leq C_0 x_0$ and
$|\lambda|<\delta_0$ then $x'\geq x_0$ provided
$|t|>\frac{C_0}{C}x+\sqrt{2x_0/C}$, $|t|<\delta_0$. Assuming $x\leq
x_0$ and taking $x_0$ sufficiently small so that
$\frac{C_0}{C}x_0+\sqrt{2x_0/C}<\delta_0$,
we thus deduce that the curve
segments $\gamma_{x,y,\lambda,\omega}|_{(-\delta_0,\delta_0)}$ are outside
the region
$x'<x_0$ for $t$ outside a (fixed!) compact subinterval of
$(-\delta_0,\delta_0)$.
{\em From now on, by $\gamma$ we mean the restriction
$\gamma_{x,y,\lambda,\omega}|_{(-\delta_0,\delta_0)}$, and we
everywhere assume that the functions we integrate along $\gamma$ are
supported in $x'\leq x_0/2$, so all integrals are on a fixed compact subinterval.}

Note that in the case of geodesics, as discussed above, $\alpha$ is a
quadratic polynomial in $\omega$; this will be of use when the ellipticity of the
boundary principal symbol is discussed.

Before we proceed, we
discuss the blowup of a space around a submanifold. Here we work
locally on say $\RR^m_w$, thus the submanifold can be taken to be
given by $w'=0$, where we write
$w=(w',w'')\in\RR^k\times\RR^{m-k}$. Then blowing up
$\RR^{m-k}=\{w'=0\}$ in $\RR^m$ amounts to introducing cylindrical
coordinates around it, i.e.\ the factor $\RR^{m-k}$ (the cylindrical
`axis', though higher dimensional) is unchanged, while on $\RR^k_{w'}$
one introduces `polar coordinates' $(|w'|,\frac{w'}{|w'|})\in
[0,\infty)\times\sphere^{k-1}$, thus one replaces $\RR^m$ by
$$
[\RR^m;\RR^{m-k}]=[0,\infty)_r\times\sphere^{k-1}_\theta\times\RR^{m-k}_{w''};
$$
altogether one has `coordinates' (the quotes are due to the spherical factor)
$$
r=|w'|,\ \theta=\frac{w'}{|w'|},\ w'',
$$
with the equalities holding outside $r=0$.
The new boundary
$$
\ff=\{0\}\times \sphere^{k-1}\times\RR^{m-k}
$$
is called
the front face.
Further one
has a blow-down map $\Phi:[\RR^m;\RR^{m-k}]\to\RR^m$ which is smooth,
namely $(r,\theta,w'')\mapsto (r\theta,w'')$, but is not invertible at
$r=0$ although it restricts to a diffeomorphism $[\RR^m;\RR^{m-k}]\setminus\ff\to\RR^m\setminus\RR^{m-k}$. We refer
to the Appendix of \cite{RBMSpec} for a concise but more detailed
description, and for further references. Note that the effect of this blow up is to distinguish directions
of approach to the submanifold being blown up, $\RR^{m-k}$; curves
$c:[0,\ep)\to\RR^m$ with $c(0)\in\RR^{m-k}$ and $c'(0)$ not in the
tangent space of $\RR^{m-k}$ lift to (i.e.\ using the diffeomorphism
property away from $\RR^{m-k}$, can be identified with) curves $\tilde
c$ in
$[\RR^m;\RR^{m-k}]$ with $\tilde c(0)\in \ff$, and two such curves
$c_j$ with $c_1(0)=c_2(0)$ satisfy $\tilde c_1(0)=\tilde c_2(0)$ if
and only if $c_1'(0)-c_2'(0)$ is tangent to $\RR^{m-k}$. (This says
that invariantly $\ff$ is the spherical normal bundle of $\RR^{m-k}$
in $\RR^m$, i.e.\ the quotient of its normal bundle minus its zero
section by dilations.)

Let $\hat X=\RR_x\times\RR^{n-1}_y$, $S\hat X=\hat X\times\RR\times\sphere^{n-2}$.
In our setting, as we show momentarily, we start with the map
\begin{equation}\label{eq:Gamma-plus-def}
\Gamma_+: S\hat X\times[0,\infty)\to [\hat X\times \hat
X;\diag],\qquad \Gamma_+(z,\nu,t)=\gamma_{z,\nu}(t),
\end{equation}
being a diffeomorphism near $S\hat X\times\{0\}$. More precisely,
$\Gamma_+$ is defined on $\tilde
U\times(-\tilde\delta_0,\tilde\delta_0)\times\sphere^{n-2}\times(-\tilde\delta_0,\tilde\delta_0)$,
and this map is a diffeomorphism onto it range when restricted to a
neighborhood of $S\hat X\times\{0\}$.
To see this, note that the
diagonal is the submanifold $z-z'=0$ of $\hat X\times\hat X$, so
nearby one can use coordinates $z-z'\in\RR^n$ (the analogue of $w'$
above) and $z\in\RR^n$ (the analogue of $w''$ above). Thus,
coordinates on $[\hat X\times \hat
X;\diag]$ are given by $z$, $|z'-z|$ and
$\frac{z'-z}{|z'-z|}$, and a simple calculation shows that at $t=0$,
one has $\frac{z'-z}{|z'-z|}=\frac{\nu}{|\nu|}$ (with the norms being
just Euclidean norms), which proves that $\Gamma_+$ as in
\eqref{eq:Gamma-plus-def} is a diffeomorphism near $S\hat X\times\{0\}$.
Similarly,
\begin{equation}\label{eq:Gamma-minus-def}
\Gamma_-: S\hat X\times(-\infty,0]\to [\hat X\times \hat
X;\diag],\qquad \Gamma_-(z,\nu,t)=\gamma_{z,\nu}(t),
\end{equation}
is a diffeomorphism near $S\hat X\times\{0\}$.

\begin{rem}
The analogous results would work
with $\tilde X$ in place of $\hat X$.
Then $S\tilde X$
is the sphere bundle of $\hat X$, i.e.\ $T\tilde X\setminus o$
quotiented out by the $\RR^+$-action. If we have a Riemannian
metric we could take this to be the unit sphere bundle with respect to
this metric, but {\em any} other choice of a transversal to the
dilation orbits in the tangent space of $\hat X$ works, such as the
unit sphere bundle with respect to another metric, or indeed (locally,
in the region of interest) the space of tangent vectors of the form
$\lambda \,\pa_x+\omega\,\pa_y$, where $\omega\in\sphere^{n-2}$,
considered above.
\end{rem}

We now reduce $\delta_0>0$ if necessary so that $\Gamma_+$ is a
diffeomorphism on $U_{x,y}\times
(-\delta_0,\delta_0)_\lambda\times\sphere^{n-2}_\omega\times
[0,\delta_0)_t$, and analogously for $\Gamma_-$; we assume this from
now on. (Note that in $\lambda_0$ we could allow an arbitrary interval
with compact closure for this particular purpose.)

Our inversion problem is now that assuming
$(If)(x,y,\lambda,\omega)=\int_{\RR}f(\gamma_{x,y,\lambda,\omega}(t))\,dt$
is known, we would like to recover $f$ from it. (Recall our
convention from above; the integral is really over
$(-\delta_0,\delta_0)$, and $f(x',y')$ vanishes for $x'\geq x_0/2$.)
It is occasionally convenient to assume
\begin{equation}\label{eq:gamma-symmetric}
\gamma_{x,y,-\lambda,-\omega}(-t)=\gamma_{x,y,\lambda,\omega}(t).
\end{equation}
Without
this symmetry assumption, we would have two curves with a given
tangent line at $(x,y)$, so having the integral of functions along
both, we would have additional information. (In other words, we could
simply drop one of these families to arrive at the present setting.)

The idea is simply
to average over the family, i.e.\ to consider for $x>0$
\begin{equation}\label{eq:A-def}
(Af)(x,y)=\int_\RR\int_{\sphere^{n-2}} (If) (x,y,\lambda,\omega)\,\tilde\chi(x,\lambda)\,d\lambda\,d\omega,
\end{equation}
where $\tilde\chi$ is supported in $|\lambda|\leq
\sqrt{2C}\sqrt{x}$. One concrete choice that achieves this
$$
\tilde\chi(x,\lambda)=x^{-1/2}\chi(\lambda/\sqrt{x}),
$$
with $\chi$ having sufficiently small support near $0$; another one is
$$
\tilde\chi(x,\lambda)=x^{-1}\chi(\lambda/x),
$$
where now any compactly supported $\chi$ works (for sufficiently small
$x$). We remark that we can allow $\chi$ to depend smoothly on
$\omega$ and $y$; over compact sets such a behavior is necessarily
uniform since there are no boundaries in these variables.

\begin{rem}
Here we need to recall that $\gamma$ and $A$ are
only locally defined, on some open set $O$ (i.e.\ $\gamma$ is defined
for $z=(x,y)\in O$ only, and only as long as its image remains in $O$). However, as we are only
interested in applying $A$ to distributions supported in $O$, and as
the ellipticity statements we show are local in nature, this
is not a problem. For instance, for $K$ a fixed subset of $O$, one may
replace $A$ by $\psi A\psi$ where $\psi\in\CI_c(O)$ is $\equiv 1$ on a
neighborhood of $K$, which is now globally well-defined, and
ellipticity statements are unaffected near $K$.
\end{rem}

For any $r$, we can write $A$ as
\begin{equation}\label{eq:A-decomp}
A=L\circ I,\ L=M_2\circ\Pi \circ M_1\circ I,
\end{equation}
where
$$
\Pi u(x,y)=\int_\RR\int_{\sphere^{n-2}} u (x,y,\lambda,\omega)\,d\lambda\,d\omega,
$$
and
$$
M_1 u(x,y,\lambda,\omega)= x^r \chi(\lambda/x) u(x,y,\lambda,\omega),\ (M_2f)(x,y)=x^{-1-r}f(x,y).
$$
Thus, $\Pi$ is a push-forward map, and thus is bounded on
$$
H^s([0,\infty)\times\RR^{n-1}\times\RR\times\sphere^{n-2})\to
H^s([0,\infty)\times\RR^{n-1})
$$
for all $s\geq 0$, i.e.\ `on $H^s$' in
brief, since such a map is bounded on $H^s$ in the
absence of boundaries, and there are continuous extension maps from
$H^s$ of a half space to $H^s$ of the whole space. On the other hand, as $\chi$ is bounded, $M_1$
is bounded on
$L^2(\{x\geq 0\})$ while its $j$th derivative is bounded by $x^{-j}$
times a constant, so $x^s M_1$ is bounded as map on $H^s(\{x\geq 0\})$
when $s\geq 0$ integer.
Thus,
$$
L:H^s([0,\infty)\times\RR^{n-1}\times\RR\times\sphere^{n-2})\to x^{-s-1}H^s([0,\infty)\times\RR^{n-1})
$$
is bounded.
Further, the X-ray transform, $I$, is itself of the form
$I=\tilde\Pi\circ\Phi^*$, where $\Phi^*$ is pull-back by the
map $(z,\nu,t)\mapsto \gamma_{z,\nu}(t)$,
$\nu=(\lambda,\omega)$, which has surjective differential
in view of the diffeomorphism property of $\Gamma_\pm$ (on the relevant
set; recall also that we are assuming that the functions we are
applying $I$ to are supported in $U$), and $\tilde\Pi$ is the
push-forward given by integration in $t$. Thus, $I$ itself is bounded
$$
I:H^s([0,\infty)\times\RR^{n-1})\to H^s([0,\infty)\times\RR^{n-1}\times\RR\times\sphere^{n-2}).
$$
Correspondingly, if we show $A$ is invertible as a map between
appropriate spaces of functions supported near $x=0$ (as
discussed in the previous section), concretely
weighted Sobolev spaces, with domain space including 
$H^s([0,\infty)\times\RR^{n-1})$ and range space including $x^{-s-1}H^s([0,\infty)\times\RR^{n-1})$,
we obtain an estimate for $f$ in
terms of $If$ when $f$ satisfies such a support condition and lies in $H^s([0,\infty)\times\RR^{n-1})$.

Note that the $A$ defined by \eqref{eq:A-def} is certainly a pseudodifferential
operator in $x>0$; moreover, its principal symbol is elliptic if
$\chi\geq 0$ with
$\chi>0$ near $0$ (this
uses $n>2$) -- while this is well-known, we check it below explicitly
in the proof of the next proposition. Our main task is to understand the uniform behavior of
$A$ to $x=0$. It turns out that while $A$ itself is not a scattering
pseudodifferential operator, its conjugates by exponential weights are:

\begin{prop}
Suppose $\chi\in\CI_c(\RR)$.
Let $\tilde\chi(x,\lambda)=x^{-1}\chi(\lambda/x)$.
The operator
$A_\digamma=x^{-1}e^{-\digamma/x}Ae^{\digamma/x}$ is in $\Psisc^{-1,0}$ for
$\digamma>0$.
\end{prop}

The main point here regarding the exponential weights is that the
Schwartz kernel of $A$ itself is well-behaved near compact subsets of
the front face, i.e.\ where $X=\frac{x'-x}{x^2}$ and $Y=\frac{y'-y}{x}$ are bounded, but is not so
well-behaved as $(X,Y)\to\infty$. However, the support conditions on
$\chi$ insure that $X\to +\infty$ on the support of the Schwartz
kernel of $A$ (with a suitable estimate), and thus the exponential
conjugation gives exponential decay of the {\em conjugated kernel} as
$(X,Y)\to\infty$, giving the conclusion of the proposition.

\begin{proof}
At first work in $x>0$, ignoring the limit $x\to 0$. Then it is
standard that $A$ is a pseudodifferential operator (the weights are
harmless then), but it is instructive to prove this in a manner that
extends seamlessly to the general case.

With $\Gamma_\pm$ as in \eqref{eq:Gamma-plus-def}-\eqref{eq:Gamma-minus-def},
for $\tilde\chi$ an arbitrary smooth function on $S\tilde X$ (not necessarily
dependent just on $x,\lambda$) the diffeomorphism property on $S\tilde
X\times [0,\delta_0)$ allows one to
rewrite, with $|d\nu|$ denoting a smooth measure on the transversal
such as $|d\lambda|\,|d\omega|$,
$$
Af(z)=\sum_{\bullet=+,-}\int f(z')\tilde\chi(\Gamma_\bullet^{-1}(z,z')) (\Gamma_\bullet^{-1})^{*}(|d\nu|\,dt)
$$
in terms of $z,z'$ as
\begin{equation}\begin{aligned}\label{eq:off-bdy-diag-ker-form}
&\int
f(z')|z'-z|^{-n+1}b\Big(z,\frac{z'-z}{|z'-z|},|z'-z|\Big)\,dz',\\
& b\Big(z,\frac{z'-z}{|z'-z|},0\Big)=\tilde\chi\Big(z,\frac{z'-z}{|z'-z|}\Big)\sigma\Big(z,\frac{z'-z}{|z'-z|}\Big),
\end{aligned}\end{equation}
where $\sigma>0$ is bounded below -- it is the change of variables
Jacobian factor. The two terms $\Gamma_\pm$ are in fact identical by
the symmetry assumption on $\gamma$, \eqref{eq:gamma-symmetric}, so we
can ignore $\Gamma_-$.
In particular, $A$ is a pseudodifferential operator with principal
symbol given by the Fourier transform of
$$
|z'-z|^{-n+1}b(z,\frac{z'-z}{|z'-z|},0)=|z'-z|^{-n+1}(\tilde\chi\sigma)(z,\frac{z'-z}{|z'-z|})
$$
in $Z=z'-z$. One can insert a cutoff $\phi$ in $|Z|$ with compact
support, identically $1$ near $0$ (considered as an even function on
$\RR$), without changing the result modulo rapid decay, i.e.\ as a
principal symbol, the result is not affected. The latter can be computed easily as
\begin{equation*}\begin{aligned}
&\int_{\RR^n} e^{-iZ\cdot\zeta}|Z|^{-n+1}(\tilde\chi\sigma)(z,\hat
Z)\phi(|Z|)\,dZ=\int_0^\infty\int_{\sphere^{n-1}}e^{-it\hat Z\cdot\zeta}(\tilde\chi\sigma)(z,\hat
Z)\phi(t)\,dt\,d\hat Z\\
&=\frac{1}{2}\int_\RR\int_{\sphere^{n-1}}e^{-it\hat Z\cdot\zeta}(\tilde\chi\sigma)(z,\hat
Z) \phi(t)\,dt\,d\hat Z=\frac{1}{2}\int_{\sphere^{n-1}}\hat\phi(\hat Z\cdot\zeta) (\tilde\chi\sigma)(z,\hat
Z)\,d\hat Z;
\end{aligned}\end{equation*}
here $\hat\phi$ is the Fourier transform of $\phi$. Fixing $\hat\zeta=\frac{\zeta}{|\zeta|}$, since $\hat\phi$
is Schwartz, if $\tilde\chi\sigma$ is supported away from the
equatorial sphere $\{\hat Z:\ \hat Z\cdot\hat\zeta=0\}$, $|\hat
Z\cdot\zeta|>c|\zeta|$ on its support for some $c>0$, and then for all $N$,
$\hat\phi(\hat Z\cdot\zeta)\leq\tilde C_N|\zeta|^{-N}$, and thus we
conclude that the integral is Schwartz and thus gives no contribution
to the principal symbol. Correspondingly (by using a partition of unity), it suffices to consider a
neighborhood of the equator and assume $\tilde\chi\sigma$ is supported
here. Then one can write $Z=(Z^\parallel,Z^\perp)$ according to the orthogonal
decomposition relative to $\hat\zeta=\frac{\zeta}{|\zeta|}$, so
$Z^\parallel=Z\cdot\hat\zeta$, similarly for $\hat Z$, and
$d\hat Z$ is of the form $a(\hat Z^\parallel)\,d\hat
Z^\parallel\,d\theta$, $\theta=\frac{\hat Z^\perp}{|\hat
  Z^\perp|}\in\sphere^{n-2}$ with $a(0)=1$ since $\hat Z^\perp=(1-|\hat
Z^\parallel|^2)^{1/2}\theta$. Thus, one has
\begin{equation*}\begin{aligned}
&\frac{1}{2}\int_\RR\int_{\sphere^{n-2}}\hat\phi(\hat Z^\parallel |\zeta|) (\tilde\chi\sigma)(z,\hat
Z^\parallel\hat\zeta+(1-|\hat
Z^\parallel|^2)^{1/2}\theta)a(\hat Z^\parallel)\,d\theta\,d\hat
Z^\parallel\\
&=\frac{1}{2|\zeta|}\int_\RR (|\zeta|\hat\phi(\hat Z^\parallel |\zeta|)) a(\hat Z^\parallel)\Big(\int_{\sphere^{n-2}}(\tilde\chi\sigma)(z,\hat
Z^\parallel\hat\zeta+(1-|\hat
Z^\parallel|^2)^{1/2}\theta)\,d\theta\Big)
\,d\hat
Z^\parallel
\end{aligned}\end{equation*}
Since $(|\zeta|\hat\phi(\hat Z^\parallel |\zeta|))\to\delta_0$ in
distributions as $|\zeta|\to\infty$, this is
$|\zeta|^{-1}\int_{\sphere^{n-2}}
  (\tilde\chi\sigma)(z,\theta)\,d\theta$ modulo terms decaying faster
  as $|\zeta|\to\infty$; indeed, one easily sees by expanding
  $\tilde\chi\sigma$
around $\hat Z^\parallel=0$ that this asymptotic holds modulo
$O(|\zeta|^{-2})$ terms. In other words, the principal symbol of $A$
at
$(z,\zeta)$ is a constant multiple of
\begin{equation}\label{eq:off-bdy-pr-symbol-form}
|\zeta|^{-1}\int (\tilde\chi\sigma)(z,\hat Z^\perp)\,dZ^\perp.
\end{equation}
In
particular, if $\tilde\chi\geq 0$, then as long as for each $(z,\zeta)$, $\zeta\neq 0$, there is $\hat Z$
perpendicular to $\zeta$ with $\tilde\chi$ non-zero at $(z,\hat Z)$,
then $A$ is an elliptic order $-1$ pseudodifferential operator, in
accordance with the results of Stefanov and Uhlmann \cite{SU2}. This is indeed
the case with our choice of $\tilde\chi$, provided $n>2$.

We now turn to the scattering behavior, i.e.\ as at least one of
$x,x'\to 0$. Note that from \eqref{eq:x-xp-comp}, on the support of
$\tilde\chi$, $x'\geq x-c_0x^2$,
for $x$ small. We in fact show below that on the support of
$\tilde\chi$, $X$ is bounded below, and $X\to+\infty$ if
$|Y|\to\infty$, and indeed $X\geq C_1|Y|^2$ for $|Y|$
sufficiently large, $C_1>0$. Here we recall from Section~\ref{sec:sc-calc} that
$$
X=\frac{x-x'}{x^2},\ Y=\frac{y-y'}{x}.
$$
With $K$ denoting the Schwartz
kernel of $A$, as
$$
x^{-1}-(x')^{-1}=\frac{x'-x}{xx'}=X\frac{x}{x'}=X/(1+xX),
$$
$A_\digamma$ has Schwartz
kernel
\begin{equation}\begin{aligned}\label{eq:K-flat-form}
K^\flat(x,y,X,Y)&=x^{-1}e^{-\digamma(x^{-1}-(x')^{-1})}K(x,y,X,Y)\\
&=x^{-1}e^{-\digamma
 X/(1+xX)}K(x,y,X,Y).
\end{aligned}\end{equation}
Taking into account the polynomial bounds on $K$ in terms of $X,Y$, and
$x'\geq x-c_0x^2$ implying that $X$ is bounded below as shown later in
the proof, further that
$X\to+\infty$ as $|Y|\to\infty$ with $X\geq C_1|Y|^2$,
exponential decay of $K^\flat$ as well as its derivatives follows
easily for $\digamma>0$. Thus, the main claim is that $K^\flat$ is smooth
for $(X,Y)$ finite, non-zero, conormal to $(X,Y)=0$.

Now, on $\Gamma_+(\supp\tilde\chi\times [0,\delta_0))$, $|x-x'|\leq
C|y-y'|$ means that locally in this region
$x,y,|y'-y|,\frac{x'-x}{|y'-y|},\frac{y'-y}{|y'-y|}$ are coordinates on
$[\tilde X\times\tilde X;\diag]$ -- indeed, this corresponds to using
the transversal $|y'-y|=1$ to dilations in
$\RR^n=\RR_{x'-x}\times\RR_{y'-y}^{n-1}$ in place of the unit sphere
$|(x'-x,y'-y)|=1$, which is indeed a transversal where $y'-y$ is large
relative to $x'-x$, i.e.\ in our region of interest.
Further, $\Gamma_+(x,y,\lambda,\omega,0)$
is, in terms of these coordinates, $(x,y,\lambda,\omega,0)$. In
general, thus,
$$
\lambda\Big(\Gamma_+^{-1}\Big(x,y,|y'-y|,\frac{x'-x}{|y'-y|},\frac{y'-y}{|y'-y|}\Big)\Big)=\Lambda\Big(x,y,|y'-y|,\frac{x'-x}{|y'-y|},\frac{y'-y}{|y'-y|}\Big),
$$
with
$$
\Lambda\Big(x,y,0,\frac{x'-x}{|y'-y|},\frac{y'-y}{|y'-y|}\Big)=\frac{x'-x}{|y'-y|},
$$
so (suppressing $\Gamma_+$ on the left hand side)
$$
\lambda=\frac{x'-x}{|y'-y|}+|y'-y|\tilde\Lambda\Big(x,y,|y'-y|,\frac{x'-x}{|y'-y|},\frac{y'-y}{|y'-y|}\Big).
$$
Now, in terms of the scattering coordinates,
$$
|y'-y|=x|Y|,\ \frac{x'-x}{|y'-y|}=\frac{xX}{|Y|},\
\frac{y'-y}{|y'-y|}=\hat Y,
$$
so (suppressing $\Gamma_+$ composed with
the scattering blow up map on the left hand side)
\begin{equation}\label{eq:lambda-over-x}
\frac{\lambda}{x}=\frac{X}{|Y|}+|Y|\tilde\Lambda\Big(x,y,x|Y|,\frac{xX}{|Y|},\hat Y\Big).
\end{equation}
Similarly,
\begin{equation*}\begin{aligned}
\omega&=\frac{y'-y}{|y'-y|}+|y'-y|\tilde
\Omega\Big(x,y,|y'-y|,\frac{x'-x}{|y'-y|},\frac{y'-y}{|y'-y|}\Big)\\
&=
\hat Y+x|Y|\tilde\Omega\Big(x,y,x|Y|,\frac{xX}{|Y|},\hat
Y\Big)
\end{aligned}\end{equation*}
and
\begin{equation}\begin{aligned}\label{eq:t-expr}
t&=|y'-y|+|y'-y|^2\tilde
T\Big(x,y,|y'-y|,\frac{x'-x}{|y'-y|},\frac{y'-y}{|y'-y|}\Big)\\
&=x|Y|+x^2|Y|^2\tilde T\Big(x,y,x|Y|,\frac{xX}{|Y|},\hat
Y\Big).
\end{aligned}\end{equation}
Thus,
\begin{equation}\begin{aligned}\label{eq:sc-density-calc}
dt\,d\lambda\,d\omega&=J\Big(x,y,\frac{X}{|Y|},|Y|,\hat Y\Big)
\,x^2|Y|^{-1} \,dX\,d|Y|\,d\hat Y\\
&=J\Big(x,y,\frac{X}{|Y|},|Y|,\hat Y\Big)\,x^2|Y|^{-n+1}\,dX\,dY
\end{aligned}\end{equation}
where the density factor $J$ is smooth and positive, and $J|_{x=0}=1$.
Also, on the blow-up
of the scattering diagonal, $\{X=0,\ Y=0\}$, in the region
$|Y|>\ep|X|$, thus on the support of $\chi$ in view of \eqref{eq:lambda-over-x},
$$
x,y,|Y|,\frac{X}{|Y|},\hat Y
$$
are valid coordinates,
with $|Y|$ being the defining function of the front face of this blow
up (i.e. of the lifted diagonal). Taking into account
the $x^{-1}$ in the definition in $\tilde\chi$, we thus deduce that
$K^\flat$ is given by
\begin{equation}\label{eq:K-flat-form}
e^{-\digamma X/(1+xX)}\chi\Big(\frac{X}{|Y|}+|Y|\tilde\Lambda(x,y,x|Y|,\frac{xX}{|Y|},\hat Y)\Big) |Y|^{-n+1}J\Big(x,y,\frac{X}{|Y|},|Y|,\hat Y\Big),
\end{equation}
so in particular it is conormal to the front face on the blow-up of the scattering
diagonal, of the form $\rho^{-n+1}b$, where $b$ is smooth up to the
front face, and without the first exponential factor it, together with
its derivatives (in $x,y,X,Y$) has polynomial
growth estimates as $(X,Y)\to\infty$, i.e.\ the derivatives satisfy
bounds $\leq C|(X,Y)|^N$ for some $C,N$ (depending on the
derivative). Decomposing $K^\flat$ into pieces supported in, say,
$|(X,Y)|<2$ and $|(X,Y)|>1$ by a partition of unity, we show in the
next paragraph that the
latter is Schwartz in $(X,Y)$ due the exponential decay of the first
factor of \eqref{eq:K-flat-form} on the support of the second
factor. On the other hand, for the former term, supported in $|(X,Y)|<2$,
calculations as in \eqref{eq:off-bdy-diag-ker-form}
in Fourier transforming this in $(X,Y)$ show that this term of $K^\flat$ is indeed
the Schwartz kernel of an element of $\Psisc^{-1,0}$, with standard
principal symbol being given by the analogue of
\eqref{eq:off-bdy-pr-symbol-form}. Here the additional information is
in the behavior at $x=0$, but given that our operator {\em is} an
element of $\Psisc^{-1,0}$, the same information can be obtained from
computing the boundary principal symbol, which we need in any case.

We use \eqref{eq:t-expr} to express $\frac{\lambda}{x}$ using
\begin{equation}\label{eq:gamma-form-near-diag}
x'=x+\lambda t+\alpha(x,y,\lambda,\omega)t^2+O(t^3),\ y'=y+\omega t+O(t^2),
\end{equation}
where the $O(t^2)$ and $O(t^3)$ terms have coefficients which are
smooth in $(x,y,\lambda,\omega)$.
Thus,
$$
X=\frac{x'-x}{x^2}=\frac{\lambda t}{x^2}+\frac{\alpha
  t^2}{x^2}+\frac{t^3}{x^2}\Upsilon(x,y,x\mu,\omega,t),
$$
with $\Upsilon$ a smooth function of its arguments,
so
\begin{equation*}\begin{aligned}
X=&\frac{\lambda(\Gamma_+^{-1})}{x}|Y|(1+x|Y|\tilde T(x,y,x|Y|,xX/|Y|,\hat Y))\\
&+\alpha(\Gamma_+^{-1}) |Y|^2(1+x|Y|\tilde T(x,y,x|Y|,xX/|Y|,\hat
Y))^2+x|Y|^3\Upsilon(\Gamma_+^{-1}),
\end{aligned}\end{equation*}
and thus
\begin{equation}\label{eq:lambda-over-x-form}
\frac{\lambda(\Gamma_+^{-1})}{x}=\frac{X-\alpha(\Gamma_+^{-1})|Y|^2}{|Y|}+O(x),
\end{equation}
where the $O(x)$ has smooth coefficients in terms of
$x,y,x|Y|,xX/|Y|,\hat Y$. Thus, for $\mu=\frac{\lambda}{x}\in (-c,c)$,
$-c|Y|<X-\alpha(\Gamma_+^{-1})|Y|^2<c|Y|$, which shows (by the positive
definiteness of $\alpha$) that $X\to+\infty$ on $\supp\tilde\chi$ if
$|Y|\to\infty$, and indeed, for $|Y|$ sufficiently large, $X>C_1|Y|^2$
for some $C_1>0$.

As already explained, this proves the proposition, since now for all
$N'$ the
exponential factor in \eqref{eq:K-flat-form} is $\leq
C'|(X,Y)|^{-N'}$ for suitable $C'$ on the support of the second
factor, so combined with the
polynomial estimates for the derivatives of the second and third
factors, it follows that $K^\flat$ is smooth in $(x,y)$, with values
in functions Schwartz in $(X,Y)$ for $(X,Y)\neq 0$, and conormal to
$(X,Y)=0$, which is exactly the characterization of the Schwartz
kernel of a scattering pseudodifferential operator.
\end{proof}

\begin{rem}
We now explain the form these arguments would take for a different
scaling chosen for $\tilde\chi$. By \eqref{eq:gamma-form-near-diag}
for $\lambda=\sqrt{x}\mu$, with $\mu$ in a compact set near $0$
(i.e.\ the first, $O(\sqrt{x})$ localization used above),
$x'=x+\sqrt{x}\mu t+\alpha(x,y,\sqrt{x}\mu,\omega,t)t^2+O(t^3)$ gives
$$
x'\leq C'(x+|y-y'|^2),
$$
indicating that $|y'-y|/x^{1/2}$ is the appropriate homogeneous
variable for analysis; using $X=\frac{\sqrt{x'}-\sqrt{x}}{\sqrt{x}}$,
$Y=\frac{y'-y}{\sqrt{x}}$, this amounts to a statement that the
analysis is well-behaved on the 0 double space of Mazzeo-Melrose \cite{Mazzeo-Melrose:Meromorphic} when
the smooth structure is given by the boundary defining function
$\sqrt{x}$. This is a somewhat complicated space with a
non-commutative normal operator at infinity; there's a reduced normal
operator after a partial Fourier transform and rescaling which is a
b-scattering (or Bessel) type pseudodifferential operator on a
half-line. This is the reason for using our sharper cutoff, which puts
us into the more amenable setting of Melrose's scattering calculus, as
described above.
\end{rem}

We now compute the boundary principal symbol of $A_\digamma$. Indeed,
this is immediate from \eqref{eq:K-flat-form} and \eqref{eq:lambda-over-x-form}
which show that at $x=0$ (i.e.\ the scattering front face) the
Schwartz kernel of $A_\digamma$ is
\begin{equation*}\begin{aligned}
e^{-\digamma
  X}|Y|^{-n+1}\chi\Big(\big(X-\alpha(0,y,0,\hat
Y)|Y|^2\big)/|Y|\Big)=\tilde K(y,X,Y).
\end{aligned}\end{equation*}
As described in Section~\ref{sec:sc-calc},
for each $y$, $\tilde K(y,.,.)$ acts as a convolution operator,
thus it becomes a multiplication operator upon Fourier transforming in $(X,Y)$,
and the desired invertibility $\Hsc^{s,r}\to \Hsc^{s+1,r}$ amounts to the Fourier
transformed kernel, $\hat K(y,.,.)$ being bounded below in absolute value by
$c\langle(\xi,\eta)\rangle^{-1}$, $c>0$ (here $(\xi,\eta)$ are the
Fourier dual variables of $(X,Y)$). Thus, we need to compute the
inverse Fourier transform of $\tilde K(y,.,.)$, and find $\digamma$ such that
the desired bound holds. Note also that if $\chi$ depends on $y$ and
$\omega$ as discussed above, we simply have
$$
\chi\Big(\big(X-\alpha(0,y,0,\hat
Y)|Y|^2\big)/|Y|,y,\hat Y\Big)
$$
in the above expression for the Schwartz kernel at the front face. We
have thus shown

\begin{lemma}
The boundary principal symbol of
$x^{-1}e^{-\digamma/x}Ae^{\digamma/x}$ is
the $(X,Y)$-Fourier transform of
$$
\tilde K(y,X,Y)=e^{-\digamma
  X}|Y|^{-n+1}\chi\Big(\big(X-\alpha(0,y,0,\hat
Y)|Y|^2\big)/|Y|,y,\hat Y\Big).
$$
\end{lemma}

In order to find a suitable $\chi$, we first make
a slightly inadmissible choice for an easier computation,
namely we take $\chi(s)=e^{-s^2/(2\nu)}$ with $\nu$ to be fixed (and
allowed to depend on $y$ and $\hat Y$), so
$\hat\chi(.)=c\sqrt{\nu}e^{-\nu|.|^2/2}$ for appropriate $c>0$. Thus,
$\chi$ does not have compact support, and an approximation argument
will be necessary.
Now, in general (for arbitrary $\chi$ which has
superexponential decay so its Fourier transform is entire), the Fourier
transform in $X$ is
\begin{equation}\label{eq:partial-FT-gen}
\cF_X \tilde K(y,\xi,Y)=|Y|^{2-n}e^{-\alpha\digamma|Y|^2}e^{-i\alpha\xi|Y|^2}\hat\chi((\xi-i\digamma)|Y|),
\end{equation}
as follows by taking into account the effect of translations,
dilations and multiplication by exponential weights on the Fourier
transform (the last two of which are closely related). Here $\alpha$
is a function of $y$ and $\hat Y$, as above. Substituting
the particular $\chi$ yields a non-zero multiple of
\begin{equation}\begin{aligned}\label{eq:partial-FT-Gauss}
&\sqrt{\nu}|Y|^{2-n}e^{-\alpha\digamma|Y|^2}e^{-i\alpha\xi|Y|^2}e^{-\nu(\xi-i\digamma)^2|Y|^2/2}\\
&=\sqrt{\nu}|Y|^{2-n}e^{-(2\alpha\digamma+\nu\xi^2-\nu\digamma^2)|Y|^2/2}e^{-i(\alpha-\digamma\nu)\xi|Y|^2}.
\end{aligned}\end{equation}
Now, the $Y$-Fourier transform of $|Y|^{2-n}$ is a homogeneous radial
(i.e.\ $\mathrm{SO}(n-1)$-invariant) function of order $-1$, so it is a
non-zero multiple of $|\eta|^{-1}$, with $\eta$ the Fourier-dual
variable of $Y$. Notice that this uses very strongly that we have
$n>2$; $n=2$ would give a delta distribution. Thus, if the $Y$-Fourier transform of
\begin{equation*}
e^{-(2\alpha\digamma+\nu\xi^2-\nu\digamma^2)|Y|^2/2}e^{-i(\alpha-\digamma\nu)|Y|^2}
\end{equation*}
is positive, then the Fourier transform of the product, which is given
by the convolution (in $\eta$) of these, is also positive, and with asymptotic
behavior given by that of $|\eta|^{-1}$ provided the Fourier transform
of the Gaussian is Schwartz.
Indeed, if one Fourier transforms $|Y|^{2-n}\psi(y, Y)$, where $\psi$
is Schwartz in the last variable, only the behavior near $Y=0$ contributes to the
asymptotics as $\eta\to\infty$, and thus using the Taylor series of
$\psi$, one obtains the asymptotic expansion of the Fourier transform
as $\eta\to\infty$ as a classical polyhomogeneous function (with the
expansion in terms of $|\eta|^{-1-j}$, $j\geq 0$ integer).

So assume now that $\alpha$ is a positive definite quadratic form
in $\hat Y$ and take $\nu=\digamma^{-1}\alpha$ (so same holds for
$\nu$, i.e.\ $\nu$ is a quadratic form in $\hat Y$).
Thus, one has $\alpha|Y|^2=Q(Y,Y)$, a quadratic form in
$Y$. Thus, writing $Q^{-1}(Y,Y)$ for the
dual quadratic form, and taking $\chi(s)=e^{-s^2/(2\digamma^{-1}Q(\hat Y,\hat
  Y))}$, we have
$$
\hat \chi(\sigma)=c(\digamma^{-1}Q(\hat Y,\hat
Y))^{1/2}e^{-\digamma^{-1}Q(\sigma\hat Y,\sigma\hat  Y)/2}.
$$
In view of \eqref{eq:partial-FT-gen}-\eqref{eq:partial-FT-Gauss},
$\cF_X \tilde K(y,\xi,Y)$ is a non-zero multiple of
\begin{equation*}\begin{aligned}
&\sqrt{\nu}|Y|^{2-n}e^{-\alpha\digamma|Y|^2}e^{-i\alpha\xi|Y|^2}e^{-\nu(\xi-i\digamma)^2|Y|^2/2}\\
&=\sqrt{\nu}|Y|^{2-n}e^{-(2\alpha\digamma+\nu\xi^2-\nu\digamma^2)|Y|^2/2}e^{-i(\alpha-\digamma\nu)\xi|Y|^2}\\
&=\digamma^{-1/2}\sqrt{\alpha}|Y|^{2-n}e^{-(\xi^2+\digamma^2)\alpha\digamma^{-1}|Y|^2/2}=\digamma^{-1/2}\sqrt{\alpha}|Y|^{2-n}e^{-\digamma^{-1}(\xi^2+\digamma^2)Q(Y,Y)/2},
\end{aligned}\end{equation*}
where the last factor is a real Gaussian
since the oscillatory factor in \eqref{eq:partial-FT-Gauss} becomes
identically $1$.
This is Schwartz in $Y$ for $\digamma>0$, and thus the Fourier
transform is a positive multiple of
$$
(\det
Q)^{-1/2}\digamma^{(n-1)/2}(\xi^2+\digamma^2)^{-(n-1)/2}e^{-\digamma Q^{-1}(\eta,\eta)/(2(\xi^2+\digamma^2))}.
$$
which satisfies the requirements from the previous paragraph (positive
Schwartz function).

One has to be a bit
careful about the joint $(\xi,\eta)$-behavior, i.e.\ when $\xi$ is
also going to infinity, and where we still need lower bounds. The
Fourier transform of $|Y|^{2-n}e^{-\digamma^{-1}\langle\xi\rangle^2Q(Y,Y)/2}$, with
$\langle\xi\rangle=(\xi^2+\digamma^2)^{1/2}$, is a constant
multiple of
\begin{equation*}\begin{aligned}
&\int
|\eta-\zeta|^{-1}\langle\xi\rangle^{-(n-1)}e^{-\digamma Q^{-1}(\zeta,\zeta)/(2\langle\xi\rangle^2)}\,d\zeta\\
&=\int
|\eta-\langle\xi\rangle\zeta'|^{-1}e^{-\digamma Q^{-1}(\zeta',\zeta')/2}\,d\zeta'\\
&=\langle\xi\rangle^{-1}\int|\eta/\langle\xi\rangle-\zeta'|^{-1}e^{-\digamma
  Q^{-1}(\zeta',\zeta')/2}\,d\zeta'
=\langle\xi\rangle^{-1}\varphi(\eta/\langle\xi\rangle)
\end{aligned}\end{equation*}
where we wrote $\zeta'=\zeta/\langle\xi\rangle$, and where $\varphi$
is an elliptic positive classical symbol of order $-1$, namely the convolution of
$|.|^{-1}$ with the Schwartz function $e^{-\digamma Q^{-1}(.,.)/2}$. This assures lower bounds
$c\langle(\xi,\eta)\rangle^{-1}$, $c>0$, i.e.\ elliptic lower
bounds. Indeed,
this is immediate when $\langle\xi\rangle>\ep|\eta|$, for
$\langle\xi\rangle^{-1}$ is equivalent to
$\langle(\xi,\eta)\rangle^{-1}$ in this region in terms of decay rates, while
$\varphi(\eta/\langle\xi\rangle)$ is a 0th order symbol in this
region.
To see what happens when $|\eta|>\ep \langle\xi\rangle$, notice that
by virtue of the classicality in fact have
$\varphi(\eta')=|\eta'|^{-1}\tilde\varphi(\langle\xi\rangle/|\eta'|,\eta'/|\eta'|)$, with
$\tilde\varphi$ smooth near $0$ in the first argument. Thus, we obtain
$$
\langle\xi\rangle^{-1}\varphi(\eta/\langle\xi\rangle)=|\eta|^{-1}\tilde\varphi(\langle\xi\rangle/|\eta|,\eta'/|\eta'|),
$$
which is a symbol of order $-1$ in $|\xi|\leq C|\eta|$, and
$|\eta|^{-1}$ is equivalent to $\langle(\xi,\eta)\rangle^{-1}$ here,
completing the proof of the ellipticity claim.

Now, if $\chi$ is not a Gaussian, but rather one has a sequence $\chi_n$
in $\CI_c(\RR)$ which converges to the Gaussian in Schwartz functions
(notice that this does not imply that the Fourier transform of
$\chi_n$ is
pointwise positive for any $n$, which is the reason we need to use the Fourier
transform of the Gaussian directly),
then the Fourier transforms converge in the appropriate spaces, which
suffices to conclude that for sufficiently large $n$, letting
$\chi=\chi_n$,
the Fourier transform $\hat K$, i.e.\ the
boundary principal symbol, still has
lower bounds $\tilde C\langle(\xi,\eta)\rangle^{-1}$, $\tilde C>0$, as
desired. We have thus proved:

\begin{lemma}
For $\digamma>0$ there exists $\chi\in\CI_c(\RR)$, $\chi\geq 0$,
$\chi(0)=1$, such that for the corresponding operator
$x^{-1}e^{-\digamma/x}Ae^{\digamma/x}$ the boundary symbol is
elliptic; indeed, this holds for all $\chi$ sufficiently close in
Schwartz space to a
specific Gaussian.
\end{lemma}

Hence, we have
$$
B=x^{-1}e^{-\digamma/x}Ae^{\digamma/x}\in\Psisc^{-1,0}
$$
elliptic both in the sense of the standard principal symbol (in the
set of interest $O$), and
the scattering principal symbol, which is at $x=0$, and in particular
the results of Section~\ref{sec:sc-calc} are applicable. Thus,
elements of the kernel of $A$ which have support in the compact subset
$K$ of $O$ is finite dimensional, and further a stability estimate
holds on a complementary subspace of this finite dimensional
subspace. Further, with $x=x_c=\tilde x+c$, as discussed at the
beginning of this section, the arguments at the end of
Section~\ref{sec:sc-calc} show that for sufficiently small $c$, this
subspace of the kernel of $A=A_c$ is actually trivial, and one has a
stability estimate in $M_c=\{x_c>0\}$ for functions supported in
$K$. Thus, for $c$ small, writing the support condition as final subscript,
$$
A=xe^{\digamma/x}Be^{-\digamma/x}:e^{\digamma/x}\Hsc^{s,r}(M_c)_K\to
xe^{\digamma/x}\Hsc^{s+1,r}(M_c)=e^{\digamma/x}\Hsc^{s+1,r+1}(M_c)
$$
satisfies estimates
$$
\|f\|_{e^{\digamma/x}\Hsc^{s,r}(M_c)_K}\leq C\|Af\|_{e^{\digamma/x}\Hsc^{s+1,r+1}(M_c)}.
$$
In particular, if one is willing to give up polynomial weights as
unimportant at the cost of losing $\delta/x$ in the exponential weight, $\delta>0$, and one uses that for $s\geq 0$, $H^s(M_c)\subset \Hsc^{s,r}(M_c)$ for
$r\leq -\frac{n+1}{2}$ while for $r\geq -\frac{n+1}{2}+2s$,
$\Hsc^{s,r}(M_c)\subset H^s(M_c)$, with continuous inclusion maps, see
\eqref{eq:reg-Sob-in-sc}-\eqref{eq:sc-Sob-in-reg}, we have
$$
\|f\|_{e^{(\digamma+\delta)/x}H^s(M_c)_K}\leq C\|Af\|_{e^{\digamma/x}H^{s+1}(M_c)}.
$$
Now, using the decomposition $A=L\circ I$ of $A$ in \eqref{eq:A-decomp}, and the boundedness
statements following it, we have for all $\digamma>0$,
$$
\|Af\|_{e^{\digamma/x}H^{s+1}(M_c)}\leq C'\|If\|_{H^{s+1}(\cM_{M_c})}
$$
when $f\in H^{s+1}(M_c)_K$.
For the convenience of the reader, we summarize all the maps for
$\digamma>0$ and with $r\in\RR$, $\delta>0$ arbitrary, in a
commutative diagram:
\begin{equation*}\begin{CD}
e^{\digamma/x}\Hsc^{s,r}(M_c)_K@>A>>e^{\digamma/x}\Hsc^{s+1,r+1}(M_c)@>G>>e^{\digamma/x}\Hsc^{s,r}(M_c)\\
@AAA @AAA @VVV\\
H^{s+1}(M_c)_K@>L\circ I>> H^{s+1}(M_c)@.e^{(\digamma+\delta)/x} H^s(M_c)
\end{CD}\end{equation*}
with all vertical arrows inclusion maps, $G$ the inverse of $A$ on the
range of $A$, and with
\begin{equation*}\begin{CD}
H^{s+1}(M_c)_K@>I>>H^{s+1}(\cM_{O_c})@>L>> \Hsc^{s+1}(M_c)
\end{CD}\end{equation*}
being the lower left composite map $L\circ I$ written out in detail.
Indeed, note that even $s=-1$ is allowed with the inclusions we
stated; factoring the first inclusion via $H^{s+1}(M_c)\to
e^{\digamma/x}\Hsc^{s+1.r}(M_c)\to e^{\digamma/x} \Hsc^{s,r}(M_c)$
proves it since $s+1\geq 0$,
while using \eqref{eq:neg-sc-Son-in-reg} (with $-s$ in place of $s$)
gives the last inclusion map.
In combination, this completes the proof of the main theorem for
$s\geq -1$ in the notation here (thus $s\geq 0$ for the notation of
the main theorem), with
$\digamma$ replaced by $\digamma+\delta$ -- as both $\digamma>0$ and
$\delta>0$ are arbitrary, this means that the original statement is proved.

\bibliographystyle{plain}
\bibliography{sm}

\end{document}